\documentclass[12pt]{amsart}

\usepackage[margin=1.1in]{geometry}
\usepackage[english]{babel}
\usepackage[utf8]{inputenc}
\usepackage{fancyhdr}
\usepackage{natbib}
\usepackage{amsmath}

\RequirePackage{amsthm,amsmath,amsfonts,amssymb}
\usepackage{graphicx}
\usepackage{hyperref}
\usepackage{graphicx}
\usepackage{placeins}
\usepackage{caption}
\usepackage{subcaption}
\numberwithin{equation}{section}

\theoremstyle{plain}
\theoremstyle{remark} 
\numberwithin{equation}{section} 

\title{Estimation of Order Restricted Location/Scale Parameters of a General Bivariate Distribution Under General Loss function: Some Unified results}

\author{Naresh Garg  and Neeraj Misra \\ {\footnotesize Department of Mathematics and Statistics\\Indian Institute of Technology Kanpur \\Kanpur-208016, Uttar Pradesh, India}}

\makeatletter
\def\@seccntformat#1{%
  \protect\textup{\protect\@secnumfont
    \ifnum\pdfstrcmp{subsection}{#1}=0 \bfseries\fi
    \csname the#1\endcsname
    \protect\@secnumpunct
  }%
}  
\makeatother

\begin{document}
\maketitle
\section*{\textbf{Abstract}}

	We consider component-wise equivariant estimation of order restricted location/scale parameters of a general bivariate distribution under quite general conditions on underlying distributions and the loss function. This paper unifies various results in the literature dealing with sufficient conditions for finding improvments over arbitrary location/scale equivariant estimators. The usefulness of these results is illustrated through various examples. A simulation study is considered to compare risk performances of various estimators under bivariate normal and independent gamma probability models. A real-life data analysis is also performed to demonstrate applicability of the derived results. 
\\~\\ \textbf{Keywords:} Improved estimator; Inadmissibility; Location equivariant estimator; Restricted MLE; Restricted Parameter Space; Scale equivariant estimator.



\section{\textbf{Introduction}}  \label{introduction}
The problems of estimating order restricted location/scale parameters (say, $\theta_1$ and $\theta_2$) of two distributions are of interest in many applications. For example, in a clinical trial, when estimating average blood pressures of two groups of Hypertension patients, one treated with a standard drug and the other with a placebo, is of interest, it can be assumed that the average blood pressure of Hypertension patients treated with the standard drug is lower than the average blood pressure of Hypertension patients treated with the placebo. Similarly, in an agricultural experiment, let $\theta_1$ be the average yield of a crop without fertilizer and $\theta_2$ be the average yield of the same crop using an effective fertilizer. Here it may be of interest to estimate $\theta_1$ and $\theta_2$ and it may be appropriate to assume that $\theta_1\leq\theta_2$ (see Kubokawa and Saleh (\citeyear{MR1370413})). An interesting application of problems of estimation of order restricted parameters in meta-analysis can be found in Taketomi et al. (\citeyear{axioms10040267}). Due to their applicability in various real-life situations, estimation of order restricted parameters has been extensively studied in the literature. For an account of developments in this area and various applications, one may refer to monographs by Barlow et al. (\citeyear{MR0326887}), Robertson et al. (\citeyear{MR961262}) and Van Eeden (\citeyear{MR2265239}). 
\vspace{2mm}

Several studies in this area deal with improving unrestricted best location/scale equivariant estimators (BLEE/BSEE) and/or unrestricted maximum likelihood estimators (MLE) by exploiting the prior information that the parameters are order restricted. However, most of these studies are centered around specific distributions, having independent marginals, and specific loss functions (see, for example, Vijayasree et al. (\citeyear{MR1345425}), Misra and Dhariyal (\citeyear{MR1326266}) and Patra and Kumar (\citeyear{Patra})). For a detailed account of other developments in this direction, one may refer to Van Eeden (\citeyear{MR2265239}). Only limited studies have been carried out for improving an arbitrary location/scale equivariant estimator and mostly these studies are centred around specific distributions, having independent marginals, and specific loss functions. Kubokawa and Saleh (\citeyear{MR1370413}) and Garg and Misra (\citeyear{garg2021componentwise}) unified some of these studies on improving best location/scale estimators when the parameters are known to satisfy an order restriction. In this paper, we attempt to further unify these studies by considering the problem of improving arbitrary location/scale equivariant estimators under general loss functions and general location/scale probability models.

The results obtained in this paper are applicable to a variety of situations that are not covered in the existing literature. Moreover, many findings already reported in the literature for specific probability distributions can be obtained as specific cases of our results. For example, estimation of order restricted means of bivariate normal distribution is not adequately studied in the literature and is covered under the umbrella of general framework of this paper. For this case, the general results derived in the paper can be used to obtain classes of estimators (also containing the restricted maximum likelihood estimators) improving upon the usual estimators under certain conditions. Here, the findings of the paper are also useful in studying the effect of the correlation coefficient on the structural form of improved estimators. We observe that, depending upon whether the correlation coefficient is smaller or larger than a ratio of standard derivations, improved estimators may underestimate or overestimate the mean (see Examples 2.1.1 and 2.2.2). When correlation coefficient equals a ratio of standard deviations, it is not possible to find improvement over the usual estimator using the results of the paper. Interestingly in this case, the BLEE and the restricted maximum likelihood estimators are the same. Similar observations may be obtained in the case of other standard bivariate distributions, including the Cheriyan bivarite gamma distribution and bivariate distributions generated through standard copulas (see Kotz et al. (\citeyear{MR1788152})).

The rest of the paper is organised as follows. In Section \ref{2}, we consider component-wise estimation of order restricted location parameters of a general bivariate location family of distributions under a general bowl-shaped loss function. Component-wise estimation of order restricted scale parameters of a general bivariate scale family of distributions under a general bowl-shaped loss function is dealt with in Section \ref{3}. Following Stein (\citeyear{MR171344}), we develop a technique for deriving sufficient conditions for inadmissibility of any location/scale equivariant estimator and obtaining their improvements. Usefulness of our results is illustrated through various examples. In Section \ref{4}, we consider a simulation study for comparing risk performances of various estimators for estimating smallest location and scale parameter under bivariate normal and independent gamma probability models, respectively. In Section \ref{5}, a real-life data analysis is considered to illustrate the applicability of derived results.

\section{\textbf{Improved Estimators for Location Parameters}} \label{2}
Let $\bold{X}=(X_1,X_2)$ be a random vector with the joint probability density function (pdf)
\begin{equation}\label{eq:2.1}
	f_{\boldsymbol{\theta}}(x_1,x_2)= 	f(x_1-\theta_1,x_2-\theta_2),\; \; \;(x_1,x_2)\in \Re^2, 
\end{equation} 
where $f(\cdot,\cdot) $ is a specified pdf on $\Re^2$ and $\boldsymbol{\theta}=(\theta_1,\theta_2)\in \Theta_0=\{(t_1,t_2)\in\Re^2:t_1 \leq t_2\}$ is the vector of unknown and restricted location parameters; here $\Re$ denotes the real line and $\Re^2=\Re\times \Re$. Generally, $\bold{X}=(X_1,X_2)$ would be a minimal-sufficient statistic based on a bivariate random sample or two independent random samples, as the case may be.  \vspace{2mm}

Consider estimation of the location parameter $\theta_i$ under a general loss function $L_i(\boldsymbol{\theta},a)=W(a-\theta_i),\;\boldsymbol{\theta}\in\Theta_0,\; a\in\mathcal{A},\;i=1,2,$     
where $\mathcal{A}=\Re$ and $W:\Re\rightarrow [0,\infty)$ is a specified non-negative function that satisfies the following conditions: \vspace*{0.5mm}

\label{(C1)}	\textbf{\boldmath(C1):} $W(0)=0$, $W(t)$ is decreasing on $(-\infty,0)$ and increasing on $(0,\infty)$; 

\label{(C2)}	\textbf{\boldmath(C2):} $W'(t)$ is non-decreasing, almost everywhere.	\vspace{2mm}

The above problem of estimating restricted location parameter $\theta_i \,(i=1,2)$ is invariant under the group of transformations $\mathcal{G}=\{g_c:\,c\in\Re\},$ where $g_c(x_1,x_2)$ $=(x_1+c,x_2+c),\; (x_1,x_2)\in\Re^2,\;c\in\Re.$ Under the group of transformations $\mathcal{G}$, any location equivariant estimator of $\theta_i$ has the form
\begin{equation}\label{eq:2.2}
	\delta_{\psi}(\bold{X})=X_i-\psi(D),
\end{equation}
for some function $\psi:\,\Re\rightarrow \Re\,,\;i=1,2,$ where $D=X_2-X_1$. Let us denote the risk function of the estimator $\delta_{\psi}(\bold{X})$, for estimating $\theta_i$, by
$R_i(\boldsymbol{\theta},\delta_{\psi})=E_{\boldsymbol{\theta}}[W(\delta_{\psi}(\bold{X})-\theta_i)],\; \, \boldsymbol{\theta}\in\Theta_0,\;i=1,2.$ Note that the risk function $R_i(\boldsymbol{\theta},\delta_{\psi}),\;i=1,2,$ depends on $\boldsymbol{\theta}\in\Theta_0$ only through $\lambda=\theta_2-\theta_1\in[0,\infty)$. Exploiting the prior information of order restriction on parameters $\theta_1$ and $\theta_2$ ($\theta_1\leq \theta_2$), we aim to derive sufficient conditions on function $\psi$ that help us to obtain estimators improving on the estimator $\delta_{\psi}$  for estimating $\theta_i,\,i=1,2$.\vspace*{2mm}

Note that, in the unrestricted case $\Theta=\Re^2$, the problem of estimating the location parameter $\theta_i$, under the loss function $L_i(\cdot,\cdot)$, is invariant under the additive group of transformations $\mathcal{G}_0=\{g_{c_1,c_2}:\,(c_1,c_2)\in\Re^2\},$ where $g_{c_1,c_2}(x_1,x_2)=(x_1+c_1,x_2+c_2),\; (x_1,x_2)\in\Re^2,\;(c_1,c_2)\in\Re^2$, and the unrestricted best location equivariant estimator (BLEE) of $\theta_i$ is $\delta_{0,i}(\bold{X})=X_i-c_{0,i}$, where $c_{0,i}$ is the unique solution of the equation $\int_{-\infty}^{\infty} \int_{-\infty}^{\infty}\, W'(s_i-c) \,f(s_1,s_2)\,ds_1\,ds_2=0,\;i=1,2.$\vspace*{2mm}

The following lemma, whose proof is straightforward, will be useful in proving the main results of the paper.
\\~\\ \textbf{Lemma 2.1.} Let $ s_0\in \Re $ be a fixed constant and let $ M_i:\Re\rightarrow [0,\infty), \; i=1,2,$ be non-negative functions such that
\begin{align*}
	M_1(s) M_2(s_0) &\geq \; (\leq)\; M_1(s_0) M_2(s),\; \forall \; s<s_0\\
	\text{and }\; \;M_1(s) M_2(s_0) &\leq\;(\geq)\; M_1(s_0) M_2(s),\; \forall \; s>s_0.
\end{align*}
Let $ M:\Re\rightarrow \Re$ be a function, such that $M(s)\leq 0,\;\forall\;s<s_0,$ and $M(s)\geq 0, \;\forall\; s>s_0$. 
Then,
$$ M_2(s_0) \int\limits_{-\infty}^{\infty} M(s) \, M_1(s) ds\leq\;(\geq)\; M_1(s_0) \int\limits_{-\infty}^{\infty} M(s) \, M_2(s) ds.$$

\subsection{\textbf{Estimation of the Smaller Location Parameter $\theta_1$}} \label{2.1}
\setcounter{equation}{0}
\renewcommand{\theequation}{2.1.\arabic{equation}}
\noindent
\vspace*{2mm}

Let $S=X_1-\theta_1$ and $\lambda=\theta_2-\theta_1$. Note that the risk function of location equivariant estimator $\delta_{\psi}(\bold{X})=X_1-\psi(D)$ of $\theta_1$ can be written as \begin{align*}
	R_1(\boldsymbol{\theta},\delta_{\psi})
	&= E_{\boldsymbol{\theta}}[W(X_1-\psi(D)-\theta_1)]\\
	&=E_{\boldsymbol{\theta}}[E_{\boldsymbol{\theta}}[W(S-\psi(D))\vert D]]\\
	&=\int_{-\infty}^{\infty} r_{1,\lambda}(\psi(t),t)\,g(t-\lambda) dt,
\end{align*}
where $r_{1,\lambda}(c,t)= E_{\lambda}[W(S-c)\vert D=t]=\int_{-\infty}^{\infty} W(s-c)f(s\vert t-\lambda)ds,$
for $f(s \vert t)=\frac{f(s,s+t)}{\int_{-\infty}^{\infty}f(s,s+t)ds}$ and $g(t)=\int_{-\infty}^{\infty} f(s,s+t)ds$. \vspace*{2mm}

The following assumption is also required for proving the main results.
\\~\\ \textbf{\boldmath(A1):} For any fixed $\lambda\geq 0$ and $t$, the equation 
\begin{equation}\label{eq:2.1.1}
	\int_{-\infty}^{\infty} \, W'(s-c)\,f(s \vert t-\lambda)\,ds=0
\end{equation}
has the unique solution $c\equiv\psi_{\lambda}(t)$, such that $f(\psi_{\lambda}(t)\vert t-\lambda)>0$. 
\\~\\For any real numbers $x$ and $y$, let $x\vee y$ and $x \wedge y$ denote their maximum and minimum, respectively.
\\~\\ \textbf{Theorem 2.1.1.} Suppose that conditions (C1), (C2) and the assumption (A1) hold. Let $\delta_{\psi}(\bold{X})=X_1-\psi(D)$ be a location equivariant estimator of $\theta_1$, where $\psi:\Re\rightarrow \Re$. Let $\underline{\psi}(t)$ and $\overline{\psi}(t)$ be functions such that $\underline{\psi}(t)<\psi_{\lambda}(t)<\overline{\psi}(t), \;\forall\;\lambda\geq 0$ and any $t$. For any fixed $t$, define $\psi^{*}(t)=(\underline{\psi}(t)\vee\psi(t))\wedge\overline{\psi}(t)=\underline{\psi}(t)\vee(\psi(t)\wedge\overline{\psi}(t))$. Then $R_1(\boldsymbol{\theta},\delta_{\psi})\geq R_1(\boldsymbol{\theta},\delta_{\psi^{*}}),$ for any $\boldsymbol{\theta}\in\Theta_0$, where $\delta_{\psi^{*}}(\bold{X})=X_1-\psi^{*}(D)$.
\begin{proof}\!\!\!. Let $A=\{t:\psi(t)<\underline{\psi}(t)\}$, $B=\{t:\underline{\psi}(t)\leq \psi(t)\leq \overline{\psi}(t)\}$ and $C=\{t:\psi(t)>\overline{\psi}(t)\}$, so that $\psi^{*}(t)=\underline{\psi}(t)\;(=\overline{\psi}(t))$, if $t\in A$ ($t\in C$) and $\psi^{*}(t)=\psi(t)$, if $t\in B$. Conditions (C1), (C2) and the assumption (A1) imply that, for any fixed $\lambda\geq 0$,
	$$r_{1,\lambda}(c,t)=\int_{-\infty}^{\infty} W(s-c)f(s\vert t-\lambda)ds $$
	is non-increasing in $c\in (-\infty,\psi_{\lambda}(t)]$, non-decreasing in $c\in [\psi_{\lambda}(t),\infty),$ with unique minimum at $c\equiv \psi_{\lambda}(t)$. Since, for any fixed $\lambda\geq 0$ and any $t$, $\underline{\psi}(t)\leq \psi_{\lambda}(t)\leq \overline{\psi}(t)$, it follows that $r_{1,\lambda}(c,t)$
	is non-increasing in $c\in (-\infty,\underline{\psi}(t)]$ and non-decreasing in $c\in [\overline{\psi}(t),\infty).$
	Consequently, $r_{1,\lambda}(\psi(t),t)\geq r_{1,\lambda}(\underline{\psi}(t),t)$, for $t\in A$, $r_{1,\lambda}(\psi(t),t)\geq r_{1,\lambda}(\overline{\psi}(t),t)$, for $t\in C$, and
	\\~\\  $R_1(\boldsymbol{\theta},\delta_{\psi})$
	\begin{align*}
		&=\int_{A}\! r_{1,\lambda}(\psi(t),t)\,g(t-\lambda) dt + \int_{B}\! r_{1,\lambda}(\psi(t),t)\,g(t-\lambda) dt + \int_{C}\! r_{1,\lambda}(\psi(t),t)\,g(t-\lambda) dt\\
		&\geq \int_{A} \!r_{1,\lambda}(\underline{\psi}(t),t)\,g(t-\lambda) dt + \int_{B}\! r_{1,\lambda}(\psi(t),t)\,g(t-\lambda) dt + \int_{C}\! r_{1,\lambda}(\overline{\psi}(t),t)\,g(t-\lambda) dt\\
		&=E_{\boldsymbol{\theta}}[W(S-\psi^{*}(D))]\\
		&=R_1(\boldsymbol{\theta},\delta_{\psi^{*}}), \; \; \boldsymbol{\theta}\in \Theta_0.
	\end{align*}
\end{proof}
In order to identify $\underline{\psi}(t)$ and $\overline{\psi}(t)$ considered in the above theorem, the following theorem may turnout to be useful in many situations.
\\~\\\textbf{Theorem 2.1.2.} Suppose that conditions (C1), (C2) and the assumption (A1) hold. If, for every fixed $\lambda\geq 0$ and $t$, $f(s,s+t-\lambda)/f(s,s+t)$ is non-decreasing (non-increasing) in $s\in \Re$, then, for every fixed $t$, $\psi_{\lambda}(t)$ is a non-decreasing (non-increasing) function of $\lambda\in [0,\infty)$.
\begin{proof}\!\!\!.
	Let us fix $\lambda_1$ and $\lambda_2$, such that $0\leq \lambda_1<\lambda_2<\infty$. Under the assumption (A1), we have, for any $t$, $f(\psi_{\lambda_i}(t)\vert t-\lambda_i)>0$ and $\int_{-\infty}^{\infty} \, W'(s-\psi_{\lambda_i}(t))\,f(s\vert t-\lambda_i)\,ds=0,\;i=1,2.$
	Suppose that $\psi_{\lambda_1}(t)\neq \psi_{\lambda_2}(t)$. Taking $s_0=\psi_{\lambda_1}(t)$, $M(s)=W'(s-\psi_{\lambda_1}(t))$, $M_1(s)=f(s\vert t-\lambda_1)$ ($M_1(s)=f(s\vert t-\lambda_2)\,$) and $M_2(s)=f(s\vert t-\lambda_2)$ ($M_2(s)=f(s\vert t-\lambda_1)\,$), and using Lemma 2.1, we conclude that
	\begin{multline*}
		0=f(\psi_{\lambda_1}(t)\vert t-\lambda_2)\,\int_{-\infty}^{\infty} \, W'(s-\psi_{\lambda_1}(t))\,f(s\vert t-\lambda_1)\,ds\, \\  
		\leq \,(\geq)\, f(\psi_{\lambda_1}(t)\vert t-\lambda_1)\, \int_{-\infty}^{\infty} \, W'(s-\psi_{\lambda_1}(t))\,f(s\vert t-\lambda_2)\,ds,
	\end{multline*}
	which, under the assumption (A1) and the assumption that $\psi_{\lambda_1}(t)\,\neq\,\psi_{\lambda_2}(t),$ further implies that $k_1(\psi_{\lambda_1}(t)\vert t)>(<)\,0$, where 
	$$k_1(c\vert y)= \int_{-\infty}^{\infty} \, W'(s-c)\,f(s\vert y-\lambda_2)\,ds,\;c\in \Re,\; y\in \Re. $$
	Since, for any $y\in \Re$, $k_1(c\vert y) \text{ is a non-increasing function of }c, \;k_1(\psi_{\lambda_1}(t)\vert t)$ $\,>\,(<)\, 0$ and $k_1(\psi_{\lambda_2}(t)\vert t)=0$, we conclude that $\psi_{\lambda_1}(t)\,<\,(>)\,\psi_{\lambda_2}(t).$ 
\end{proof}	 

Define, for any fixed $t$, \begin{equation}\label{eq:2.1.2}
	\underline{\psi}(t) =\inf_{\lambda\geq 0} \psi_{\lambda}(t)
	\qquad \text{and} \qquad
	\overline{\psi}(t)=\sup_{\lambda\geq 0} \psi_{\lambda}(t).
\end{equation}

The proof of the following corollary is contained in the proof of Theorem 2.1.1.
\\~\\	\textbf{Corollary 2.1.1.} Suppose that conditions (C1), (C2) and the assumption (A1) hold and let $\delta_{\psi}(\bold{X})=X_1-\psi(D)$ be a location equivariant estimator of $\theta_1$, where $\psi:\Re\rightarrow \Re$. Let $\psi_{0,1}:\Re\rightarrow \Re$ be such that $\psi(t)\leq \psi_{0,1}(t) \leq \underline{\psi}(t)$, whenever $\psi(t)\leq \underline{\psi}(t)$, and $\overline{\psi}(t)\leq \psi_{0,1}(t)\leq \psi(t)$, whenever $\overline{\psi}(t)\leq \psi(t)$, where $\underline{\psi}(\cdot)$ and $\overline{\psi}(\cdot)$ are defined by \eqref{eq:2.1.2}. Also let $\psi_{0,1}(t)=\psi(t)$, whenever $\underline{\psi}(t)\leq \psi(t) \leq \overline{\psi}(t)$. Then, 
$R_1(\boldsymbol{\theta},\delta_{\psi_{0,1}})$ $\leq R_1(\boldsymbol{\theta},\delta_{\psi}),\; \forall \; \boldsymbol{\theta}\in \Theta_0,$
where $\delta_{\psi_{0,1}}(\bold{X})=X_1-\psi_{0,1}(D)$.
\\~\\\textbf{Remark 2.1.1.} Under the unrestricted parameter space $\Theta=\Re^2$, the unrestricted best location equivariant estimator (BLEE) of $\theta_1$ is
$\delta_{0,1}(\bold{X})$ $=X_1-c_{0,1},$ where  $c_{0,1}$ is the unique solution of the equation
$\int\limits_{-\infty}^{\infty} \int\limits_{-\infty}^{\infty}\! W'(s_1-c)\,f(s_1,s_2)ds_1 ds_2=0$. Suppose that conditions (C1), (C2) and the assumption (A1) hold $\text{and }  P_{\boldsymbol{\theta}}[c_{0,1}\in [\underline{\psi}(D),\,\overline{\psi}(D)]\,]<1,$ for some  $\boldsymbol{\theta}\in \Theta_0$. Then, using Theorem 2.1.1, we conclude that the unrestricted BLEE $\delta_{0,1}(\bold{X})$ is inadmissible  for estimating $\theta_1$ and is dominated by the estimator $\delta_{\psi_{0,1}^{*}}(\bold{X})=X_1-\psi_{0,1}^{*}(D),$ where $\psi_{0,1}^{*}(D)$ $=(\underline{\psi}(D)\vee c_{0,1})\wedge\overline{\psi}(D)=\underline{\psi}(D)\vee(c_{0,1}\wedge\overline{\psi}(D))$. Kubokawa and Saleh (\citeyear{MR1370413}) and Garg and Misra (\citeyear{garg2021componentwise}) have provided conditions under which the condition $P_{\boldsymbol{\theta}}[c_{0,1}\in [\underline{\psi}(D),\,\overline{\psi}(D)]\,]<1,$ for some $\boldsymbol{\theta}\in \Theta_0$, holds.
\\~\\Now we will provide some applications of Theorem 2.1.1 and Corollary 2.1.1.
\\~\\ \textbf{Example 2.1.1.} Let $\bold{X}=(X_1,X_2)$ follow a bivariate normal distribution with joint pdf given by (\ref{eq:2.1}), where, for known positive real numbers $\sigma_1$ and $\sigma_2$ and known $\rho \in (-1,1),$ the joint pdf of $(X_1-\theta_1,X_2-\theta_2)$ is
$$ f(z_1,z_2) =\frac{1}{2 \pi \sigma_1 \sigma_2 \sqrt{1-\rho^2}} e^{-\frac{1}{2(1-\rho^2)}\left[\frac{z_1^2}{\sigma_1^2}-2 \rho \, \frac{z_1 z_2}{\sigma_1 \sigma_2}+\frac{z_2^2}{\sigma_2^2}\right]},\; \; \; \bold{z}=(z_1,z_2)\in \Re^2.$$ 
Consider estimation of $\theta_1$ under the squared error loss function (i.e., $W(t)=t^2,\; t\in\Re$). Let $\psi_{\lambda}(t)$ be as defined in the assumption (A1). It can be verified that\vspace*{1mm}

\begin{small}
	$$\psi_{\lambda}(t)=\frac{\int_{-\infty}^{\infty}\,s\,f(s,s+t-\lambda)ds}{\int_{-\infty}^{\infty}\,f(s,s+t-\lambda)ds}=\frac{\sigma_1 (\rho \sigma_2-\sigma_1)}{\sigma_1^2+\sigma_2^2-2\rho \sigma_1 \sigma_2}\,(t-\lambda),\;\;\; \lambda\geq 0,\;t\in\Re .$$
\end{small}

Let $\delta_{\psi}(\bold{X})=X_1-\psi(D)$ be a location equivariant estimator of $\theta_1$ and let $\psi^{*}$ be as defined in Theorem 2.1.1. If $\rho \sigma_2<\sigma_1$, then $\psi_0(t)\leq \psi_{\lambda}(t)<\infty$ and $\psi^{*}(t)=\psi(t)\vee\psi_0(t)$. If $\rho \sigma_2>\sigma_1$, then $-\infty<\psi_{\lambda}(t)\leq \psi_0(t)$ and $\psi^{*}(t)=\psi(t)\wedge\psi_0(t)$. Using Theorem 2.1.1, with
$$ \psi^{*}(t)=\begin{cases}
	\psi(t)\vee\psi_0(t),&\text{if}\;\; \rho\sigma_2<\sigma_1\\
	\psi(t)\wedge\psi_0(t),&\text{if}\;\; \rho\sigma_2>\sigma_1
\end{cases},$$
it follows that if $P_{\boldsymbol{\theta}}[\psi(D)\neq \psi^{*}(D)]>0$, for some $\boldsymbol{\theta}\in\Theta_0$, then the estimator $\delta_{\psi}$ is inadmissible for estimating $\theta_1$ and is dominated by $\delta_{\psi^{*}}(D)=X_1-\psi^{*}(D)$. One can apply these findings to obtain improvements over mixed estimators (or isotonic regression) defined in Kumar and Sharma (\citeyear{MR981031}). 

The unrestricted BLEE of $\theta_1$ is $\delta_{0,1}(\bold{X})=X_1$ (i.e., $\psi_{0,1}(t)=0,\;t\in\Re$). Then, the BLEE $\delta_{0,1}$ is improved on by 
\begin{equation}\label{eq:2.1.3}
	\delta_{1,RMLE}(\bold{X})= \min\Big\{X_1, \frac{\sigma_2(\sigma_2-\rho\sigma_1)X_1+\sigma_1(\sigma_1-\rho\sigma_2)X_2}{\sigma_1^2+\sigma_2^2-2\rho \sigma_1 \sigma_2}\Big\},
\end{equation}
when $\rho \sigma_2<\sigma_1$, and the BLEE is improved on by
\begin{equation}\label{eq:2.1.4}
	\delta_{1,RMLE}(\bold{X})= \max\Big\{X_1, \frac{\sigma_2(\sigma_2-\rho\sigma_1)X_1+\sigma_1(\sigma_1-\rho\sigma_2)X_2}{\sigma_1^2+\sigma_2^2-2\rho \sigma_1 \sigma_2}\Big\},
\end{equation}
when $\rho \sigma_2>\sigma_1$. It is easy to verify that $\delta_{1,RMLE}$ is the restricted maximum likelihood estimator of $\theta_1$ under the restricted parameter space $\Theta_0$ (see Patra and Kumar (\citeyear{Patra})). Note that, when $\rho \sigma_2=\sigma_1$, we are not able to get improvements over the BLEE using our results. Interestingly, in this case, the BLEE and the restricted maximum likelihood estimator are the same.
\\~\\ \textbf{Example 2.1.2.} Let $X_1$ and $X_2$ be dependent random variables with the joint pdf
(\ref{eq:2.1}), where the joint pdf of $(X_1-\theta_1,X_2-\theta_2)$ is
$f(z_1,z_2)= 2 z_1 z_2 \,e^{-z_1}e^{-z_2} ,\text{ if } 0< z_1<z_2<\infty;\,=0,\text{ otherwise}$. \vspace*{1.5mm}

For estimation of $\theta_1$, consider the LINEX loss function
$L_1(\boldsymbol{\theta},a)\!=\!W(a-\theta_1),\;a\in\mathcal{A}=\Re,\; \boldsymbol{\theta}\in\Theta_0,$ where $W(t)=e^t-t-1, \;t\in\Re$. Here $W(t)$ satisfies conditions (C1) and (C2). We have, for $t\in [\lambda,\infty)$ and $\lambda\geq 0$,
\begin{align*}
	\psi_{\lambda}(t)
	=\ln\left(\frac{\int_{-\infty}^{\infty}e^{s}\,f(s,s+t-\lambda)ds}{\int_{-\infty}^{\infty}\,f(s,s+t-\lambda)ds}\right)
	= \ln\left(\frac{4(2+t-\lambda)}{1+t-\lambda}\right).	\end{align*}

Under (\ref{eq:2.1.2}), we have $\underline{\psi}(t)=\!	\psi_{0}(t)=\!	\ln\left(\frac{4(2+t)}{1+t}\right)
\;\text{and}\;
\overline{\psi}(t)\! 
=\lim\limits_{\lambda\to t}	\psi_{\lambda}(t)=\ln(8).$\vspace*{2mm}

Using Corollary 2.1.1, we conclude that any location equivariant estimator $\delta_{\psi}(\bold{X})=X_1-\psi(D)$, with $P_{\boldsymbol{\theta}}[\psi(D)\neq \psi_{0}(D)\vee (\psi(D)\wedge\ln(8))]>0,$ for some $\boldsymbol{\theta}\in\Theta_0$, is inadmissible for estimating $\theta_1$ and is improved on by $\delta_{\psi^{*}}(\bold{X})=X_1-\psi^{*}(D)$, where $	\psi^{*}(D)=\max\big\{\ln\left(\frac{4(2+D)}{1+D}\right),\min\{\psi(D),\ln(8)\}\big\}.$
\\~\\ The unrestricted BLEE of $\theta_1$ is $\delta_{0,1}(\bold{X})=X_1-c_{0,1}$, where $c_{0,1}=\ln(E(e^{X_1-\theta_1}))=\ln(6)$. From the above discussion it follows that the unrestricted BLEE $\delta_{0,1}(\bold{X})=X_1-\ln(6)$ is inadmissible for estimating $\theta_1$ and is improved on by
$\delta_{\psi_{0,1}^{*}}(\bold{X})\!\!=\!\!X_1\!-\!\max\left\{\ln\left(\frac{4(2+D)}{1+D}\right),\ln(6)\right\}.$
\\~\\ \textbf{Example 2.1.3.} Let $X_1$ and $X_2$ be independent random variables with the joint pdf $f(x_1-\theta_1,x_2-\theta_2)= \frac{1}{\sigma_1 \sigma_2}\, e^{-\frac{x_1-\theta_1}{\sigma_1}}e^{-\frac{x_2-\theta_2}{\sigma_2}},\;x_1>\theta_1,\;x_2>\theta_2,\;(\theta_1,\theta_2)\in\Theta_0$, where $\sigma_1$ and $\sigma_2$ are known positive constants. Here, for every fixed $t\in \Re$ and $\lambda\geq 0$, $f(s,s+t-\lambda)/f(s,s+t)$ is non-decreasing in $s\in[0,\infty)$. \vspace*{1.5mm}

Consider estimation of $\theta_1$ under the squared error loss function ($W(t)=t^2,\;t\in\Re$).
We have, for $t\in\Re$, 
\begin{align*}
	\psi_{\lambda}(t)
	&=\frac{\int_{-\infty}^{\infty}\,s\,f(s,s+t-\lambda)ds}{\int_{-\infty}^{\infty}\,f(s,s+t)ds}
	= \max\{0,-t+\lambda\}+\frac{\sigma_1 \sigma_2}{\sigma_1+\sigma_2},\;\lambda\geq 0.
\end{align*}
Clearly, under the notations of Corollary 2.1.1,
$\underline{\psi}(t)=\!	\psi_{0}(t)= \max\{0,-t\}+\frac{\sigma_1 \sigma_2}{\sigma_1+\sigma_2}
\text{ and }
\overline{\psi}(t)\! 
=\lim\limits_{\lambda\to \infty}	\psi_{\lambda}(t)=\infty,\;t\in\Re.$ Let $\delta_{\psi}(\bold{X})=X_1-\psi(D)$ be an equivariant estimator of $\theta_1$ with $P_{\boldsymbol{\theta}}[\psi(D)\neq \psi_0(D)]>0$, for some $\boldsymbol{\theta}\in\Theta_0$. Using Corollary 2.1.1, it follows that the equivariant estimator $\delta_{\psi}(\bold{X})=X_1-\psi(D)$ is inadmissible for estimating and is dominated by $\delta_{\psi^{*}}(\bold{X})=X_1-\psi^{*}(D)$, where
$\psi^{*}(D)=\max\{\underline{\psi}(D),\psi(D)\}.$
\vspace*{2mm}

The restricted MLE of $\theta_1$ is $\delta_{1,RMLE}(\bold{X})=\min\{X_1,X_2\}$. Using Corollary 2.1.1, it follow that the restricted MLE $\delta_{\psi_{1,RMLE}}$ is inadmissible for estimating $\theta_1$ and the estimator
$\delta_{1,RMLE}^{*}(\bold{X})=\min\{X_1,X_2\}-\frac{\sigma_1 \sigma_2}{\sigma_1+\sigma_2}$
dominates it.\vspace*{2mm}

Under the squared error loss function, the unrestricted BLEE of $\theta_1$ is $\delta_{0,1}(\bold{X})=X_1-\sigma_1$.
It follows that the unrestricted BLEE $\delta_{0,1}(\bold{X})$ is dominated by $\delta_{\psi_{0,1}^{*}}(\bold{X})=X_1-\psi_{0,1}^{*}(D)$, where
$\psi_{0,1}^{*}(t)=\max\{\psi_0(t),\sigma_1\}=\max\big\{\frac{\sigma_1\sigma_2}{\sigma_1+\sigma_2}-t,\sigma_1\big\},\;\forall\; t\in\Re.$ Note that $ \delta_{\psi_{0,1}^{*}}(\bold{X})=\min\big\{X_2-\frac{\sigma_1\sigma_2}{\sigma_1+\sigma_2},X_1-\sigma_1\big\}.$\vspace*{1.5mm}

The above findings, which are easy consequences of Corollary 2.1.1, are also reported in Vijayasree et al. (\citeyear{MR1345425}). Garren (\citeyear{MR1802627}) has shown the universal dominance of the estimator $\delta_{\psi_{0,1}^{*}}(\bold{X})$ over the unrestricted BLEE $\delta_{0,1}(\bold{X})=X_1-\sigma_1$.

\subsection{\textbf{Estimation of the Larger Location Parameter $\theta_2$}} \label{2.2}
\setcounter{equation}{0}
\renewcommand{\theequation}{2.2.\arabic{equation}}
\noindent
\vspace*{2mm}

Consider estimation of the larger location parameter $\theta_2$ under the general loss function $L_2(\boldsymbol{\theta},a)=W(a-\theta_2),\;\boldsymbol{\theta}\in\Theta_0,\; a\in\mathcal{A}= \Re$, when it is known apriori that $\boldsymbol{\theta}\in\Theta_0$. Assume that $W(\cdot)$ satisfies conditions (C1) and (C2). The form of any location equivariant estimator of $\theta_2$ is $\delta_{\psi}(\bold{X})=X_2-\psi(D),$ for some function $\psi:\,\Re\rightarrow \Re$, where $D=X_2-X_1$.\vspace*{2mm} 

Let $S=X_2-\theta_2$ and $\lambda=\theta_2-\theta_1$. Note that the risk function of location equivariant estimator $\delta_{\psi}(\bold{X})=X_2-\psi(D)$ of $\theta_2$ can be written as \begin{align*}
	R_2(\boldsymbol{\theta},\delta_{\psi})
	=\int_{-\infty}^{\infty} r_{2,\lambda}(\psi(t),t)\,g(t-\lambda) dt,
\end{align*}
where $r_{2,\lambda}(c,t)= E_{\lambda}[W(S-c)\vert D=t]=\int_{-\infty}^{\infty} W(s-c)f(s\vert t-\lambda)ds,$
for $f(s \vert t)=\frac{f(s-t,s)}{\int_{-\infty}^{\infty}f(s-t,s)ds}$ and $g(t)=\int_{-\infty}^{\infty} f(s-t,s)ds$. \vspace*{2mm}

The following assumption is required for proving the main results.
\\~\\ \textbf{\boldmath(A2):} For any fixed $\lambda\geq 0$ and $t$, the equation 
\begin{equation}\label{eq:2.2.1}
	\int_{-\infty}^{\infty} \, W'(s-c)\,f(s \vert t-\lambda)\,ds=0
\end{equation}
has the unique solution $c\equiv\psi_{\lambda}(t)$, such that $f(\psi_{\lambda}(t)\vert t-\lambda)>0$. \vspace*{2mm}

Now we will state the main results of this subsection. The proofs of the following theorems and corollary are similar to those of Theorems 2.1.1-2.1.2 and Corollary 2.1.1, respectively, and hence skipped.
\\~\\ \textbf{Theorem 2.2.1.} Suppose that conditions (C1), (C2) and the assumption (A2) hold. Let $\delta_{\psi}(\bold{X})=X_2-\psi(D)$ be a location equivariant estimator of $\theta_2$, where $\psi:\Re\rightarrow \Re$. Let $\underline{\psi}(t)$ and $\overline{\psi}(t)$ be functions such that $\underline{\psi}(t)<\psi_{\lambda}(t)<\overline{\psi}(t), \;\forall\;\lambda\geq 0$ and $t$. For any fixed $t$, let $\psi^{*}(t)=(\underline{\psi}(t)\vee\psi(t))\wedge\overline{\psi}(t)=\underline{\psi}(t)\vee(\psi(t)\wedge\overline{\psi}(t))$. Then, $R_2(\boldsymbol{\theta},\delta_{\psi})\geq R_2(\boldsymbol{\theta},\delta_{\psi^{*}}),$ for any $\boldsymbol{\theta}\in\Theta_0$, where $\delta_{\psi^{*}}(\bold{X})=X_2-\psi^{*}(D)$.
\\~\\\textbf{Theorem 2.2.2.} Suppose that conditions (C1), (C2) and the assumption (A2) hold. If, for every fixed $\lambda\geq 0$ and $t$, $f(s-t+\lambda,s)/f(s-t,s)$ is non-increasing (non-decreasing) in $s\in \Re$, then  $\psi_{\lambda}(t)$ is a non-increasing (non-decreasing) function of $\lambda\in [0,\infty)$.
\\~\\Define, for any fixed $t$, \begin{equation}\label{eq:2.2.2}
	\underline{\psi}(t) =\inf_{\lambda\geq 0} \psi_{\lambda}(t)
	\qquad \text{and} \qquad
	\overline{\psi}(t)=\sup_{\lambda\geq 0} \psi_{\lambda}(t).
\end{equation}
\textbf{Corollary 2.2.1.} Suppose that conditions (C1), (C2) and the assumption (A2) hold and let $\delta_{\psi}(\bold{X})=X_2-\psi(D)$ be a location equivariant estimator of $\theta_2$, where $\psi:\Re\rightarrow \Re$. Let $\psi_{0,2}:\Re\rightarrow \Re$ be such that $\psi(t)\leq \psi_{0,2}(t) \leq \underline{\psi}(t)$, whenever $\psi(t)\leq \underline{\psi}(t)$, and $\overline{\psi}(t)\leq \psi_{0,2}(t)\leq \psi(t)$, whenever $\overline{\psi}(t)\leq \psi(t)$, where $\underline{\psi}(\cdot)$ and $\overline{\psi}(\cdot)$ are defined by (\ref{eq:2.2.2}). Also let $\psi_{0,2}(t)=\psi(t)$, whenever $\underline{\psi}(t)\leq \psi(t) \leq \overline{\psi}(t)$. Then, 
$R_2(\boldsymbol{\theta},\delta_{\psi_{0,2}})$ $\leq R_2(\boldsymbol{\theta},\delta_{\psi}),\; \forall \; \boldsymbol{\theta}\in \Theta_0,$
where $\delta_{\psi_{0,2}}(\bold{X})=X_2-\psi_{0,2}(D)$.
\\~\\ As in Section 2.1, we will now apply Theorems 2.2.1-2.2.2 and Corollary 2.2.1 to estimation of the larger location parameter $\theta_2$ in probability models considered in Examples 2.1.1-2.1.3.
\\~\\ \textbf{Example 2.2.1.} Let $\bold{X}=(X_1,X_2)$ have the bivariate normal distribution as described in Example 2.1.1. Consider estimation of $\theta_2$ under the squared error loss function $L_2(\boldsymbol{\theta},a)=(a-\theta_2)^2,\;\boldsymbol{\theta}\in\Theta_0,\; a\in \mathcal{A}=\Re$. For $t\in \Re$, $\psi_{\lambda}(t)=\frac{\sigma_2 (\sigma_2-\rho \sigma_1)(t-\lambda)}{\sigma_1^2+\sigma_2^2-2\rho \sigma_1\sigma_2},\;\lambda\geq 0.$ The unrestricted BLEE of $\theta_2$ is $\delta_{0,2}(\bold{X})=X_2$ (i.e., $\psi_{0,2}(t)=0,\;t\in\Re$). Then, the BLEE $\delta_{0,2}$ is improved on by 
$$\delta_{2,RMLE}(\bold{X})= \max\Big\{X_2, \frac{\sigma_2(\sigma_2-\rho\sigma_1)X_1+\sigma_1(\sigma_1-\rho\sigma_2)X_2}{\sigma_1^2+\sigma_2^2-2\rho \sigma_1 \sigma_2}\Big\},$$ 
when $\rho \sigma_1<\sigma_2$, and the BLEE is improved on by
$$\delta_{2,RMLE}(\bold{X})= \min\Big\{X_2, \frac{\sigma_2(\sigma_2-\rho\sigma_1)X_1+\sigma_1(\sigma_1-\rho\sigma_2)X_2}{\sigma_1^2+\sigma_2^2-2\rho \sigma_1 \sigma_2}\Big\},$$ 
when $\rho \sigma_1>\sigma_2$. It is easy to verify that $\delta_{2,RMLE}$ is the restricted maximum likelihood estimator of $\theta_2$ under the restricted parameter space $\Theta_0$ (see Patra and Kumar (\citeyear{Patra})). Note that, when $\rho \sigma_1=\sigma_2$, we are not able to get improvements over the BLEE using our results. Interestingly, in this case, the BLEE and the restricted maximum likelihood estimator are the same.
\\~\\ \textbf{Example 2.2.2.} Let $X_1$ and $X_2$ be dependent random variables as described in Example 2.1.2. Consider estimation of $\theta_2$ under the LINEX loss function $L_2(\boldsymbol{\theta},a)=W(a-\theta_2),\;\boldsymbol{\theta}\in\Theta_0,\; a\in \mathcal{A}=\Re$, where 
$W(t)=e^{t}-t-1,\; t\in\Re$. We have, for $t\in[\lambda,\infty)$ and $\lambda\geq 0$,
\begin{align*}
	\psi_{\lambda}(t)
	=\ln\left(\frac{\int_{-\infty}^{\infty}e^{s}\,f(t-\lambda-s,s)\,ds}{\int_{-\infty}^{\infty}\,f(t-\lambda-s,s)\,ds}\right)
	=t-\lambda +\ln\left(\frac{4(2+t-\lambda)}{1+t-\lambda}\right).
\end{align*}
Under (\ref{eq:2.2.2}), we have
$\underline{\psi}(t)=\ln(8)$ and $
\overline{\psi}(t)=t+\ln\left(\frac{4(2+t)}{1+t}\right).$
\vspace*{2mm}

Now, using Corollary 2.2.1, we conclude that any location equivariant estimator $\delta_{\psi}(\bold{X})=X_2-\psi(D)$, with $P_{\boldsymbol{\theta}}[\psi(D)\neq \ln(8)\vee (\psi(D)\wedge\overline{\psi}(D))]>0,$ for some $\boldsymbol{\theta}\in\Theta_0$, is inadmissible for estimating $\theta_2$ and is improved on by $\delta_{\psi^{*}}(\bold{X})=X_2-\psi^{*}(D)$, where $	\psi^{*}(D)=\max\big\{\ln(8),\min\big\{\psi_2(D),D+\ln\left(\frac{4(2+D)}{1+D}\right)\big\}\big\}.$
Note that, here, the unrestricted BLEE of $\theta_2$ does not exists as $E[e^{Z_2}]=\infty$.
\\~\\\textbf{Example 2.2.3.} Let $X_1$ and $X_2$ be independent exponential random variables as considered in Example 2.1.3. Consider estimation of $\theta_2$ under the squared error loss function $L_2(\boldsymbol{\theta},a)=(a-\theta_2)^2,\;\boldsymbol{\theta}\in\Theta_0,\; a\in \mathcal{A}=\Re$. For $\lambda\geq 0$ and $t\in\Re$, we have
$$\psi_{\lambda}(t)= \max\{0,t-\lambda\}+\frac{\sigma_1 \sigma_2}{\sigma_1+\sigma_2}.$$
Clearly,
$\underline{\psi}(t)=\frac{\sigma_1 \sigma_2}{\sigma_1+\sigma_2}
\;\;	\text{and}\;\;
\overline{\psi}(t)= \frac{\sigma_1 \sigma_2}{\sigma_1+\sigma_2}+\max\{t,0\}.$
\vspace*{1.5mm}

Using Corollary 2.2.1, it follows that any location equivariant estimator $\delta_{\psi_2}(\bold{X})=X_2-\psi_2(D)$, with $P_{\boldsymbol{\theta}}[\psi(D)\neq \ln(8)\vee (\psi(D)\wedge\overline{\psi}(D))]>0,$ for some $\boldsymbol{\theta}\in\Theta_0$, is inadmissible for estimating $\theta_2$ and the estimator $\delta_{\psi_2^{*}}(\bold{X})=X_2-\psi_2^{*}(D)$ dominates it, where $\psi_2^{*}(D)=\max\big\{\frac{\sigma_1\sigma_2}{\sigma_1+\sigma_2},\min\big\{\psi_2(D),\max\{0,D\}+\frac{\sigma_1\sigma_2}{\sigma_1+\sigma_2}\big\}\big\}.$\vspace*{1.5mm}

The restricted MLE of $\theta_2$ is $\delta_{\psi_2,M}(\bold{X})=X_2$ and it is dominated by the unrestricted BLEE $\delta_{0,2}(\bold{X})=X_2-\sigma_2$. By Corollary 2.2.1, the unrestricted BLEE $\delta_{0,2}(\bold{X})=X_2-\sigma_2$ is inadmissible for estimating $\theta_2$ and an improved estimator is
\begin{align*}
	\delta_{\psi_{0,2}^{*}}(\bold{X})
	=\begin{cases} X_2-\frac{\sigma_1\sigma_2}{\sigma_1+\sigma_2},& X_2< X_1\\ X_1-\frac{\sigma_1\sigma_2}{\sigma_1+\sigma_2},& X_1\leq X_2< X_1+\frac{\sigma_2^2}{\sigma_1+\sigma_2}\\ X_2-\sigma_2,& X_2\geq X_1+\frac{\sigma_2^2}{\sigma_1+\sigma_2} \end{cases}.
\end{align*}
The above results are also derived by Vijayasree et al. (\citeyear{MR1345425}).

\section{\textbf{Improved Estimators for Scale Parameters}} \label{3}
\setcounter{equation}{0}
\renewcommand{\theequation}{3.\arabic{equation}}
Let $\bold{X}=(X_1,X_2)$ be a random vector with the joint pdf
\begin{equation}\label{eq:3.1}
	f_{\boldsymbol{\theta}}(x_1,x_2)= \frac{1}{\theta_1 \theta_2}	f\left(\frac{x_1}{\theta_1},\frac{x_2}{\theta_2}\right),\; \; \;(x_1,x_2)\in \Re^2, 
\end{equation} 
where $f(\cdot,\cdot) $ is a specified pdf, $\boldsymbol{\theta}=(\theta_1,\theta_2)\in \Theta_0=\{(t_1,t_2)\in\Re_{++}^2:t_1 \leq t_2\}$ is the vector of unknown restricted scale parameters and $\Re_{++}^2=(0,\infty)\times(0,\infty)$. Throughout, we assume that the distributional support of $\bold{X}=(X_1,X_2)$ is a subset of $\Re_{++}^2$.  Generally, $\bold{X}=(X_1,X_2)$ would be a minimal-sufficient statistic based on a bivariate random sample or two independent random samples.  \vspace*{1.5mm}

For estimation of $\theta_i$, under the restricted parameter space $\Theta_0$, consider the loss function $L_i(\boldsymbol{\theta},a)=W(\!\frac{a}{\;\theta_i}\!),\; \boldsymbol{\theta}\in\Theta_0,\; \; a\in \mathcal{A}=\Re_{++},\;i=1,2,$     
where $\Re_{++}=(0,\infty)$ and $W:[0,\infty)\rightarrow [0,\infty)$ is a function satisfying the following two conditions:\vspace*{0.5mm}

\textbf{\boldmath(C1):} $W(1)=0$, $W(t)$ is decreasing on $(0,1]$ and increasing on $[1,\infty)$;

\textbf{\boldmath(C2):} $W'(t)$ is non-decreasing on $\Re_{++}$.\vspace*{1.5mm}

The problem of estimating $\theta_i$, under the restricted parameter space $\Theta_0$ and the loss function $L_i,i\!=\!1,2,$ is invariant under the group of transformations $\mathcal{G}\!=\!\{g_b:\,b\in(0,\infty)\},$ where $g_b(x_1,x_2)\!=\!(b\,x_1,b\,x_2),\; (x_1,x_2)\in\Re^2,\;b\in(0,\infty)$. For $T=\frac{X_2}{X_1}$, any scale equivariant estimator of $\theta_i$ has the form $\delta_{\psi}(\bold{X})\!=\!\psi(T) X_i,$ for some function $\psi\!:\!\Re_{++}\!\rightarrow \Re_{++},\,i=1,2.$ The risk function of the equivariant estimator $\delta_{\psi}$ of $\theta_i$ is given by $R_i(\boldsymbol{\theta},\delta_{\psi})=E_{\boldsymbol{\theta}}\left[W\left(\frac{\delta_{\psi}(\bold{X})}{\theta_i}\right)\right]\!, \, \boldsymbol{\theta}\in\Theta_0,$
and it depends on $\boldsymbol{\theta}\in\Theta_0$ only through $\lambda=\frac{\theta_2}{\theta_1}\in[1,\infty)$. \vspace*{1.5mm}

In the unrestricted case (i.e., when $\Theta=\Re_{++}^2$), the problem of estimating $\theta_i$, under the loss function $L_i(\boldsymbol{\theta},a)=W(\!\frac{a}{\;\theta_i}\!),\; \boldsymbol{\theta}\in\Theta,\; \; a\in \mathcal{A}=\Re_{++},$ is invariant under the multiplicative group of transformations $\mathcal{G}_0\!=\!\{g_{b_1,b_2}\!:\!\,(b_1,b_2)\!\in\Re_{++}^2\},$ where $g_{b_1,b_2}(x_1,x_2)=(b_1x_1,b_2x_2),\; (x_1,x_2)\in\Re^2,\;(b_1,b_2)\in\Re_{++}^2$, and the best scale equivariant estimator of $\theta_i$ is $\delta_{0,i}(\bold{X})=c_{0,i}X_i$, where $c_{0,i}$ is the unique solution of the equation $\int_{0}^{\infty} \int_{0}^{\infty}\,s_i \,W'(cs_i) \,f(s_1,s_2)ds_1 ds_2=0,i\!=\!1,2.$ \vspace*{2mm}

In the following two subsections, we consider component-wise estimation of order restricted scale parameters $\theta_1$ and $\theta_2$ and derive sufficient conditions that ensure inadmissibility of scale equivariant estimators. For equivariant estimators satisfying the sufficient conditions, we also provide dominating estimators. Applications of main results are illustrated through various examples dealing with specific probability models. \vspace{2mm}

\subsection{\textbf{Estimation of The Smaller Scale Parameter $\theta_1$}} \label{3.1}
\setcounter{equation}{0}
\renewcommand{\theequation}{3.1.\arabic{equation}}
\noindent
\vspace*{2mm}

Consider a scale equivariant estimator $\delta_{\psi}(\bold{X})=X_1\psi(T)$, for some function $\psi:\,\Re_{++}\rightarrow \Re_{++}$. Define $S=\frac{X_1}{\theta_1}$, $T=\frac{X_2}{X_1}$ and $\lambda=\frac{\theta_2}{\theta_1}\;(\lambda\geq 1)$. The risk function of $\delta_{\psi}(\bold{X})=X_1\psi(T)$ can be written as
\begin{align*}
	R_1(\boldsymbol{\theta},\delta_{\psi})=& E_{\boldsymbol{\theta}}\left[W\left(\frac{X_1\psi(T)}{\theta_1}\right)\right]\\
	=&E_{\boldsymbol{\theta}}[E_{\boldsymbol{\theta}}[W(\psi(T)S)\vert T]],
\end{align*}
where
$$E_{\boldsymbol{\theta}}[W(S\psi(T))\vert T=t]=\int_{0}^{\infty} \, W(\psi(t)s)\,f\!\left(s\vert \frac{t}{\lambda}\right)\,ds,$$

for $f\!\left(s\vert t\right)=\frac{sf\left(s,st\right)}{\int_{0}^{\infty}sf\left(s,st\right)ds },\;s\in\Re_{++}$. \vspace*{1.5mm}

We will need the following assumption for proving the results of this section.
\\~\\ \textbf{\boldmath(A3):} For any fixed $\lambda\geq 1$ and $t$, the equation 
\begin{equation}
	\int_{0}^{\infty} \, s\,W'(cs)\,f\!\left(s\vert \frac{t}{\lambda}\right)\,ds =0
\end{equation}\label{eq:3.1.1}
has the unique solution $c\!\equiv\!\psi_{\lambda}(t)$, such that $f\!\left(\psi_{\lambda}(t)\vert \frac{t}{\lambda}\right)>0$.
\vspace*{2mm}

Firstly, we state the following result that provides estimators dominating any scale equivariant estimator of $\theta_1$ under some conditions. The proof of the theorem is similar to that of Theorem 2.1.1 and hence it is being omitted.
\\~\\ \textbf{Theorem 3.1.1.} Suppose that conditions (C1), (C2) and the assumption (A3) hold. Let $\delta_{\psi}(\bold{X})=X_1\psi(T)$ be a scale equivariant estimator of $\theta_1$, where $\psi:\Re_{++}\rightarrow \Re_{++}$. Let $\underline{\psi}(t)$ and $\overline{\psi}(t)$ be functions such that $\underline{\psi}(t)<\psi_{\lambda}(t)<\overline{\psi}(t),\,\forall\;\lambda\geq 1$ and $t$. For any fixed $t$, define $\psi^{*}(t)=\underline{\psi}(t)\vee(\psi(t)\wedge\overline{\psi}(t))=(\underline{\psi}(t)\vee\psi(t))\wedge\overline{\psi}(t)$. Then, $	R_1(\boldsymbol{\theta},\delta_{\psi})\geq 	R_1(\boldsymbol{\theta},\delta_{\psi^{*}}),$ for any $\boldsymbol{\theta}\in\Theta_0$, where $\delta_{\psi^{*}}(\bold{X})=X_1\psi^{*}(T)$.

\vspace*{1.5mm}

The following theorem deals with monotonicity of the function $\psi_{\lambda}(t)$, with respect to $\lambda\geq 1$, and can be used to determine $\underline{\psi}(t)$ and $\overline{\psi}(t)$ considered in Theorem 3.1.1. The proof of the theorem, being similar to that of Theorem 2.1.2, is omitted.
\\~\\\textbf{Theorem 3.1.2.} Suppose that conditions (C1), (C2) and the assumption (A3) hold. If, for every fixed $t$ and $\lambda\geq 1$, $f\left(s,\frac{st}{\lambda}\right)/f(s,st)$ is non-decreasing (non-increasing) in $s$, then  $\psi_{\lambda}(t)$ is a non-increasing (non-decreasing) function of $\lambda\in[1,\infty)$.
\\~\\  For any fixed $t$, define
\begin{equation}\label{eq:3.1.2}
	\underline{\psi}(t)=\inf\limits_{\lambda\geq 1} \psi_{\lambda}(t)\quad
	\text{ and }\quad
	\overline{\psi}(t)=\sup\limits_{\lambda\geq 1} \psi_{\lambda}(t).
\end{equation}
\textbf{Corollary 3.1.1.} Suppose that conditions (C1) and (C2) hold and let $\delta_{\psi}(\bold{X})=X_1\psi(T)$ be a scale equivariant estimator of $\theta_1$, where $\psi:\Re_{++}\rightarrow \Re_{++}$. Let $\psi_{0,1}:\Re_{++}\rightarrow \Re_{++}$ be such that $\psi(t)\leq \psi_{0,1}(t) \leq \underline{\psi}(t)$, whenever $\psi(t)\leq \underline{\psi}(t)$, and $\overline{\psi}(t)\leq \psi_{0,1}(t)\leq \psi(t)$, whenever $\overline{\psi}(t)\leq \psi(t)$, where $\underline{\psi}(t)$ and $\overline{\psi}(t)$ are defined by (\ref{eq:3.1.2}). Also let $\psi_{0,1}(t)=\psi(t)$, whenever $\underline{\psi}(t)\leq \psi(t) \leq \overline{\psi}(t)$. Then, $	R_1(\boldsymbol{\theta},\delta_{\psi})\geq 	R_1(\boldsymbol{\theta},\delta_{\psi_{0,1}}),\;\forall \;\underline{\theta}\in\Theta_0$,
where $\delta_{\psi_{0,1}}(\bold{X})=X_1\psi_{0,1}(T)$. \vspace*{2mm}

Now we will consider some applications of Theorem 3.1.1 and Corollary 3.1.1 to specific probability models and specific loss functions.
\\~\\\textbf{Example 3.1.1.} Let $X_1$ and $X_2$ be dependent random variables with joint pdf (\ref{eq:3.1}),
where $\boldsymbol{\theta}=(\theta_1,\theta_2)\in\Theta_0$ is a vector of unknown restricted scale parameters and $f(z_1,z_2)= 
e^{-\max\{z_1,z_2\}} (1-e^{-\min\{z_1,z_2\}}),\text{ if } (z_1,z_2)\in\Re_{++}^2;
=0,  \text{ otherwise}.  $
This bivariate distribution is a special case of Cheriyan and Ramabhadran's bivariate gamma distribution (see Kotz et al. (\citeyear{MR1788152})). Here random variable $\frac{X_1}{\theta_1}$ and $\frac{X_2}{\theta_2}$ are identically  distributed with the common gamma pdf $f(z)=z\,e^{-z},\;\text{if } z>0;=0,\text{ otherwise}$.\vspace*{2mm}

Consider estimation of the scale parameter $\theta_1$, under the squared error loss function (i.e., $W(t)=(t-1)^2,\;t>0$). Let $\psi_{\lambda}(t)$ be as defined in assumption (A3). From (\ref{eq:3.1.1}), for any $\lambda\geq 1$, we have 
\begin{align*}\psi_{\lambda}(t)= \frac{\int_{-\infty}^{\infty} s\,f\!\left(s\vert \frac{t}{\lambda}\right)ds}{\int_{-\infty}^{\infty}s^2 f\!\left(s\vert \frac{t}{\lambda}\right)ds}
	&= \begin{cases} \frac{1}{3}\frac{1-\frac{\lambda}{(t+\lambda)^3}}{1-\frac{\lambda}{(t+\lambda)^4}}, \qquad \;\; 0<t<\lambda \\~\\ \frac{1}{3}\frac{\frac{1}{t^3}-\frac{1}{(t+\lambda)^3}}{\frac{\lambda}{t^4}-\frac{\lambda}{(t+\lambda)^4}},\qquad t\geq \lambda \end{cases}. 
\end{align*}
It is easy to verify that assumptions of Theorem 3.1.2 are satisfied and, therefore, $\underline{\psi}(t)=\inf\limits_{\lambda\geq 1} \psi_{\lambda}(t)=\lim\limits_{\lambda \to \infty} \psi_{\lambda}(t)=\frac{1}{4}$ and
\\$ 
\overline{\psi}(t)=\sup\limits_{\lambda\geq 1} \psi_{\lambda}(t)=\psi_1(t)= \begin{cases} \frac{1}{3}\frac{1-\frac{1}{(t+1)^3}}{1-\frac{1}{(t+1)^4}}, \qquad \;\; 0<t<1 \\~\\\frac{1}{3} \frac{\frac{1}{t^3}-\frac{1}{(t+1)^3}}{\frac{1}{t^4}-\frac{1}{(t+1)^4}},\qquad t\geq 1 \end{cases}.$\vspace*{2mm}

By an application of Theorem 3.1.1, any equivariant estimator $\delta_{\psi}(\bold{X})=X_1\psi(T)$, with $P_{\boldsymbol{\theta}}[\frac{1}{4}<\psi(T)<\overline{\psi}(T)]>0$, for some $\boldsymbol{\theta}\in\Theta_0$, is inadmissible for estimating $\theta_1$. In this case a dominating estimator is given by  $\delta_{\psi^{*}}(\bold{X})=X_1\psi^{*}(T)$, where $\psi^{*}(T)=\max\big\{\frac{1}{4},\min\big\{\psi(T),\overline{\psi}(T)\big\}\big\}.$
\vspace*{1mm}

The unrestricted BSEE of $\theta_1$ is $\delta_{0,1}(\bold{X})=\frac{X_1}{3} $. Define, $\psi_{0,1}^{*}(T)=\max\{\underline{\psi}(T),\min\{\frac{1}{3},\overline{\psi}(T)\}\}=\min\big\{\frac{1}{3},\overline{\psi}(T)\big\}.$ Using Theorem 3.1.1, it follows that the unrestricted BSEE $\delta_{0,1}(\bold{X})=\frac{X_1}{3}$ is inadmissible for estimating $\theta_1$ and is dominated by the estimator $\delta_{\psi_{0,1}^{*}}(\bold{X})=X_1\psi_{0,1}^{*}(T)$.
\\~\\ \textbf{Example 3.1.2.} Let $X_1$ and $X_2$ be independent random variables with joint pdf (\ref{eq:3.1}),
where $\boldsymbol{\theta}=(\theta_1,\theta_2)\in \Theta_0$ and $f\left(z_1,z_2\right)= \alpha_1 \alpha_2 z_1^{\alpha_1-1}z_2^{\alpha_2-1},\text{ if } 0<z_1 <1,\,0<z_2 <1; =0, \text{ otherwise}$, for known positive constants $\alpha_1$ and $\alpha_2$.
\vspace*{1.5mm}

Consider estimation of $\theta_1$ under the squared error loss function $L_1(\boldsymbol{\theta},a)=(\frac{a}{\theta_1}-1)^2,\, \boldsymbol{\theta}\in\Theta_0,\, a\in \mathcal{A}=\Re_{++}.$ Let $\psi_{\lambda}(t)$ be as defined in assumption (A3), so that, for $\lambda\geq 1$ and $t\in \Re_{++}$,
$$\psi_{\lambda}(t)=\frac{\int_{-\infty}^{\infty} s\,f\!\left(s\vert \frac{t}{\lambda}\right)ds}{\int_{-\infty}^{\infty}s^2 f\!\left(s\vert \frac{t}{\lambda}\right)ds}
=\frac{\alpha_1+\alpha_2+2}{\alpha_1+\alpha_2+1}\max\bigg\{1,\frac{t}{\lambda}\bigg\}.$$

Then, for $t\in [0,\infty)$, $\underline{\psi}(t)=\inf\limits_{\lambda\geq 1} \psi_{\lambda}(t)=\frac{\alpha_1+\alpha_2+2}{\alpha_1+\alpha_2+1}$ and
$\overline{\psi}(t)=\sup\limits_{\lambda\geq 1} \psi_{\lambda}(t)=\frac{\alpha_1+\alpha_2+2}{\alpha_1+\alpha_2+1}\max\{1,t\}.$
\\~\\For any equivariant estimator $\delta_{\psi}(\bold{X})=X_1\psi(T)$, define,
$$\psi^{*}(T)=\max\!\big\{\!\frac{\alpha_1+\alpha_2+2}{\alpha_1+\alpha_2+1},\min\big\{\psi(T),\frac{\alpha_1+\alpha_2+2}{\alpha_1+\alpha_2+1}\max\{1,T\}\!\big\}\!\big\}.$$ Using Theorem 3.1.1, it follows that the equivariant estimator $\delta_{\psi}(\bold{X})=X_1\psi(T) $ is dominated by $\delta_{\psi^{*}}(\bold{X})=X_1\psi^{*}(T)$, provided $P_{\boldsymbol{\theta}}[\psi(T)\neq\psi^{*}(T)]>0$, for some $\boldsymbol{\theta}\in\Theta_0$.\vspace*{1.5mm}

The unrestricted BSEE is $\delta_{0,1}(\bold{X})=\frac{\alpha_1+2}{\alpha_1+1} X_1=X_1\psi_{0,1}(T)$ (say), $\text{where }\psi_{0,1}(t)=\frac{\alpha_1+2}{\alpha_1+1},\;t\in\Re_{++}$. Note that
\\$P_{\boldsymbol{\theta}}[\underline{\psi}(T)\leq \psi_{0,1}(T)\leq \overline{\psi}(T)]=P_{\boldsymbol{\theta}}\left[T\geq \frac{(\alpha_1+2)(\alpha_1+\alpha_2+1)}{(\alpha_1+1)(\alpha_1+\alpha_2+2)}\right]<1,\;\forall\;\boldsymbol{\theta}\in\Theta_0.$\vspace*{1.5mm}

Consequently, using Theorem 3.1.1, we conclude that the unrestricted BSEE $\delta_{0,1}(\bold{X})=\frac{\alpha_1+2}{\alpha_1+1} X_1$ is inadmissible for estimating $\theta_1$ and is dominated by $\delta_{\psi_{0,1}^{*}}(\bold{X})=\psi_{0,1}^{*}(T)X_1, $
where

$$ \psi_{0,1}^{*}(T)=\max\Bigg\{\frac{\alpha_1+\alpha_2+2}{\alpha_1+\alpha_2+1},\min\Big\{\frac{\alpha_1+2}{\alpha_1+1},\frac{\alpha_1+\alpha_2+2}{\alpha_1+\alpha_2+1}\max\{1,T\}\Big\}\Bigg\}$$
\newpage
		$$=\begin{cases} \frac{\alpha_1+\alpha_2+2}{\alpha_1+\alpha_2+1},&T<1\\~\\
			\frac{\alpha_1+\alpha_2+2}{\alpha_1+\alpha_2+1}T, &1\leq T\leq \frac{(\alpha_1+2)(\alpha_1+\alpha_2+1)}{(\alpha_1+1)(\alpha_1+\alpha_2+2)}\\~\\
			\frac{\alpha_1+2}{\alpha_1+1}, &T>\frac{(\alpha_1+2)(\alpha_1+\alpha_2+1)}{(\alpha_1+1)(\alpha_1+\alpha_2+2)} \end{cases}.\qquad\qquad\qquad$$	
Misra and Dhariyal (\citeyear{MR1326266}) derived similar results.
\\~\\\textbf{Example 3.1.3.} Let $X_1$ and $X_2$ be independent gamma random variables with joint pdf (\ref{eq:3.1}), where, $f(z_1,z_2)=\frac{z_1^{\alpha_1-1}z_2^{\alpha_2-1}e^{-z_1}e^{-z_2}}{\Gamma(\alpha_1)\Gamma(\alpha_2)},\; (z_1,z_2)\in\Re_{++}^2$, for known positive constants $\alpha_1$ and $\alpha_2.$\vspace*{1.5mm}

Consider estimation of the smaller scale parameter $\theta_1$, under the squared error loss function (i.e., $W(t)=(t-1)^2,\;t\in\Re_{++}).$
\\~\\ From (\ref{eq:3.1.1}), we have 
$$\psi_{\lambda}(t)=\frac{\int_{-\infty}^{\infty} s\,f\!\left(s\vert \frac{t}{\lambda}\right)ds}{\int_{-\infty}^{\infty}s^2 f\!\left(s\vert \frac{t}{\lambda}\right)ds}=\frac{\left(1+\frac{t}{\lambda}\right)}{\alpha_1+\alpha_2+1}, \; t\in (0,\infty),\; \lambda\geq 1.$$
Consequently, for any $t>0$, $\underline{\psi}(t)=\inf\limits_{\lambda\geq 1} \psi_{\lambda}(t)=\frac{1}{\alpha_1+\alpha_2+1}
\;\;\text{and}\;\; 
\overline{\psi}(t)=\sup\limits_{\lambda\geq 1} \psi_{\lambda}(t)=\frac{1+t}{\alpha_1+\alpha_2+1}.$\vspace*{1.5mm}

By an application of Theorem 3.1.1, we conclude that any equivariant estimator $\delta_{\psi}(\bold{X})=X_1\psi(T) $ with $P_{\boldsymbol{\theta}}[\underline{\psi}(T)<\psi(T)<\overline{\psi}(T)]<1$, for some $\boldsymbol{\theta}\in\Theta_0$, is inadmissible for estimating $\theta_1$ and the estimator $\delta_{\psi^{*}}(\bold{X})=X_1\psi^{*}(T)$ improves upon it, where $\psi^{*}(T)=\max\big\{\frac{1}{\alpha_1+\alpha_2+1},\min\big\{\psi(T),\frac{1+T}{\alpha_1+\alpha_2+1}\big\}\big\}.$
\\~\\$\text{Consequently, the restricted MLE } \text{ }\delta_{\psi_{1,RMLE}}\!(\bold{X})\!=\!\min\big\{\!\frac{X_1}{\alpha_1},\frac{X_1+X_2}{\alpha_1+\alpha_2}\!\big\}\!=X_1\min\big\{\frac{1}{\alpha_1}.\frac{1+T}{\alpha_1+\alpha_2}\}$ \\(see Kaur and Singh (\citeyear{MR1128873}) and Vijayasree et al. (\citeyear{MR1345425})) is inadmissible for estimating $\theta_1$ and an improved estimator is $\delta_{\psi_{1,RMLE}^{*}}(\bold{X})=X_1\max\big\{\frac{1}{\alpha_1+\alpha_2+1},\min\big\{\frac{1}{\alpha_1},\frac{1+t}{\alpha_1+\alpha_2+1}\}\}= \min\big\{\frac{X_1}{\alpha_1},\frac{X_1+X_2}{\alpha_1+\alpha_2+1}\big\}$.\vspace*{2mm}

The unrestricted BSEE of $\theta_1$ is $\delta_{0,1}(\bold{X})=\frac{X_1}{\alpha_1+1}=X_1\psi_{0,1}(T)$ (say), where $\psi_{0,1}(t)=\frac{1}{\alpha_1+1},\;t>0$. Define, 
\\$\psi_{0,1}^{*}(T)=\max\{\underline{\psi}(T),\min\{\psi_{0,1}(T),\overline{\psi}(T)\}\}=\min\Big\{\frac{1}{\alpha_1+1},\frac{1+T}{\alpha_1+\alpha_2+1}\Big\}.$
\\~\\Since, $P_{\boldsymbol{\theta}}[\underline{\psi}(T)\leq \psi_{0,1}(T)\leq \overline{\psi}(T)]=P_{\boldsymbol{\theta}}\left[T\geq \frac{\alpha_2}{\alpha_1+1}\right]<1,\;\forall\;\boldsymbol{\theta}\in\Theta_0,$
using Theorem 3.1.1, it follows that the unrestricted BSEE $\delta_{0,1}(\bold{X})=\frac{X_1}{\alpha_1+1}$ is inadmissible for estimating $\theta_1$ and it is dominated by the estimator $\delta_{\psi_{0,1}^{*}}(\bold{X})=X_1\psi_{0,1}^{*}(T)=\min\big\{\frac{X_1}{\alpha_1+1},\frac{X_1+X_2}{\alpha_1+\alpha_2+1}\}$. 
\\Similar results are also obtained by Vijayasree et al. (\citeyear{MR1345425}).

\subsection{\textbf{Estimation of The Larger Scale Parameter $\theta_2$}} \label{3.2}
\setcounter{equation}{0}
\renewcommand{\theequation}{3.2.\arabic{equation}}
\noindent
\vspace*{2mm}

For estimation of $\theta_2$ under the prior information $\boldsymbol{\theta}\in\Theta_0$ and loss function $L_2(\cdot,\cdot)$, any scale equivariant estimator of $\theta_2$ is of the form
$\delta_{\psi}(\bold{X})=X_2\psi(T),$ for some function $\psi:\,\Re_{++}\rightarrow \Re_{++}$, where $T=\frac{X_2}{X_1}$.
\\~\\Let $S=\frac{X_2}{\theta_2}$. The risk function of $\delta_{\psi}(\bold{X})=X_2\psi(T)$ is
\begin{align*}
	R_2(\boldsymbol{\theta},\delta_{\psi})
	= E_{\boldsymbol{\theta}}\left[W\left(\frac{X_2\psi(T)}{\theta_2}\right)\right]
	=E_{\boldsymbol{\theta}}[E_{\boldsymbol{\theta}}[W(\psi(T)S)\vert T]], \;\;\; \boldsymbol{\theta}\in\Theta_0,
\end{align*}
where
$$E_{\boldsymbol{\theta}}[W(S\psi(T))\vert T=t]=\int_{0}^{\infty} \, W(\psi(t)s)\,f\!\left(s\vert \frac{t}{\lambda}\right)\,ds,$$

for $f\!\left(s\vert t\right)=\frac{sf\left(\frac{s}{t},s\right)}{\int_{0}^{\infty}sf\left(\frac{s}{t},s\right)ds },\;s\in\Re_{++}$. \vspace*{1.5mm}

We will need the following assumption for the results of this section.
\\~\\ \textbf{\boldmath(A4):} For any fixed $\lambda\geq 1$ and $t$, the equation 
\begin{equation}\label{eq:3.2.1}
	\int_{0}^{\infty} \, s\,W'(cs)\,f\!\left(s\vert \frac{t}{\lambda}\right)\,ds =0
\end{equation}
has the unique solution $c\!\equiv\!\psi_{\lambda}(t)$, such that $f\!\left(\psi_{\lambda}(t)\vert \frac{t}{\lambda}\right)>0$.
\vspace*{2mm}

The proofs of following theorems and corollary are on lines of Theorem 2.1.1, Theorem 2.1.2 and Corollary 2.1.1, and, therefore, are omitted. 
\\~\\ \textbf{Theorem 3.2.1.} Suppose that conditions (C1), (C2) and the assumption (A4) hold. Let $\delta_{\psi}(\bold{X})=X_2\psi(T)$ be a scale equivariant estimator of $\theta_2$, where $\psi:\Re_{++}\rightarrow \Re_{++}$. Let $\underline{\psi}(t)$ and $\overline{\psi}(t)$ be functions such that $\underline{\psi}(t)<\psi_{\lambda}(t)<\overline{\psi}(t),\,\forall\;\lambda\geq 1$ and $t$. For any fixed $t$, define $\psi^{*}(t)=\underline{\psi}(t)\vee(\psi(t)\wedge\overline{\psi}(t))=(\underline{\psi}(t)\vee\psi(t))\wedge\overline{\psi}(t)$. Then, $	R_2(\boldsymbol{\theta},\delta_{\psi})\geq 	R_2(\boldsymbol{\theta},\delta_{\psi^{*}}),$ for any $\boldsymbol{\theta}\in\Theta_0$, where $\delta_{\psi^{*}}(\bold{X})=X_2\psi^{*}(T)$.
\\~\\\textbf{Theorem 3.2.2.} Suppose that conditions (C1), (C2) and the assumption (A4) hold. If, for every fixed $t$ and $\lambda\geq 1$, $f\left(\frac{s\lambda}{t},s\right)/f\left(\frac{s}{t},s\right)$ is non-decreasing (non-increasing) in $s$, then  $\psi_{\lambda}(t)$ is a non-increasing (non-decreasing) function of $\lambda\in[1,\infty)$.
\\~\\  For any fixed $t$, define
\begin{equation}\label{eq:3.2.2}
	\underline{\psi}(t)=\inf\limits_{\lambda\geq 1} \psi_{\lambda}(t)\quad
	\text{ and } \quad
	\overline{\psi}(t)=\sup\limits_{\lambda\geq 1} \psi_{\lambda}(t).
\end{equation}
\textbf{Corollary 3.2.1.} Suppose that conditions (C1) and (C2) hold and let $\delta_{\psi}(\bold{X})=X_2\psi(T)$ be a scale equivariant estimator of $\theta_2$, where $\psi:\Re_{++}\rightarrow \Re_{++}$. Let $\psi_{0,2}:\Re_{++}\rightarrow \Re_{++}$ be such that $\psi(t)\leq \psi_{0,2}(t) \leq \underline{\psi}(t)$, whenever $\psi(t)\leq \underline{\psi}(t)$, and $\overline{\psi}(t)\leq \psi_{0,2}(t)\leq \psi(t)$, whenever $\overline{\psi}(t)\leq \psi(t)$, where $\underline{\psi}(t)$ and $\overline{\psi}(t)$ are defined by (\ref{eq:3.2.2}). Also let $\psi_{0,2}(t)=\psi(t)$, whenever $\underline{\psi}(t)\leq \psi(t) \leq \overline{\psi}(t)$. Then, $	R_2(\boldsymbol{\theta},\delta_{\psi})\geq 	R_2(\boldsymbol{\theta},\delta_{\psi_{0,2}}),\;\forall\;\underline{\theta}\in\Theta_0$,
where $\delta_{\psi_{0,2}}(\bold{X})=X_2\psi_{0,2}(T)$. \vspace*{2mm}

Now we will consider applications of Theorems 3.2.1-3.2.2 and Corollary 3.2.1 to estimation of the larger scale parameter $\theta_2$ in probability models considered in Examples 3.1.1-3.1.3.
\\~\\\textbf{Example 3.2.1.} Let $(X_1,X_2)$ be a vector of dependent random variables as defined in Example 3.1.1. 
We have, from (\ref{eq:3.2.1}),
	$$\psi_{\lambda}(t)= \begin{cases} \frac{\lambda}{3t} \frac{1-\left(\frac{\lambda}{\lambda+t}\right)^3}{1-\left(\frac{\lambda}{\lambda+t}\right)^4}, & 0<t\leq \lambda \\~\\ \frac{1}{3} \frac{1-\left(\frac{t}{\lambda+t}\right)^3}{1-\left(\frac{t}{\lambda+t}\right)^4},& t> \lambda \end{cases}, \;  \lambda\geq 1.$$

Here, for every fixed $t\in(0,\infty)$, $\psi_{\lambda}(t)$ an increasing function of $\lambda\in[0,\infty)$. Consequently, 
	\begin{align*}
		\underline{\psi}(t) =\begin{cases} \frac{1}{3t} \frac{1-\left(\frac{1}{\lambda+t}\right)^3}{1-\left(\frac{1}{1+t}\right)^4}, & 0<t\leq 1 \\~\\ \frac{1}{3} \frac{1-\left(\frac{1}{1+t}\right)^3}{1-\left(\frac{1}{1+t}\right)^4},& t> 1 \end{cases}
		\;\;\;\;\;\;\text{and}\; \;\;\;\;\;
		\overline{\psi}(t)=\infty,\; t>0.	
	\end{align*}
By Theorem 3.2.1, any equivariant estimator $\delta_{\psi_2}(\bold{X})=\psi_2(T)X_2$, with $P_{\boldsymbol{\theta}}[\psi_2(T)\geq \underline{\psi}(T) ]<1,$ for some $\boldsymbol{\theta}\in\Theta_0$, is inadmissible for estimating $\theta_2$ and is dominated by $\delta_{\psi_2^{*}}(\bold{X})=\psi_2^{*}(T)X_2$, where $\psi_{2}^{*}(T)
=\max\big\{\underline{\psi}(T),\,\psi_2(T)\big\}.$\vspace*{2mm}

Consequently, the unrestricted BSEE $\delta_{0,2}(\bold{X})=\frac{X_2}{3} $ is inadmissible for estimating $\theta_2$ and is improved upon by the estimator $\delta_{\psi_{0,2}^{*}}(\bold{X})=\psi_{0,2}^{*}(T)X_2$, where $\psi_{0,2}^{*}(T)\!=\!\max\big\{\underline{\psi}(T),\frac{1}{3}\!\big\}.$
\\~\\ \textbf{Example 3.2.2.} Let $X_1$ and $X_2$ be independent random variables as described in Example 3.1.2. Consider estimation of $\theta_2$ under the squared error loss function
$ L_2(\boldsymbol{\theta},a)=(\frac{a}{\theta_2}-1)^2\, , \; \boldsymbol{\theta}\in \Theta_0,\;a\in (0,\infty).$ For $\lambda\geq 1$ and $t\in [0,\infty)$, we have
$$\psi_{\lambda}(t)
=\frac{\alpha_1+\alpha_2+2}{\alpha_1+\alpha_2+1}\max\bigg\{1,\frac{\lambda}{t}\bigg\}.$$
From (\ref{eq:3.2.2}), for $t\in [0,\infty)$, $\underline{\psi}(t)=\frac{\alpha_1+\alpha_2+2}{\alpha_1+\alpha_2+1}$ $\max\big\{1,\frac{1}{t}\big\}
\;\;\text{and}\;\;
\overline{\psi}(t)=\infty.$\vspace*{1.5mm}

Using Theorem 3.2.1, we conclude that any equivariant estimator $\delta_{\psi_2}(\bold{X})=\psi_2(T) X_2$, with $P_{\boldsymbol{\theta}}[\underline{\psi}(T)\geq \psi_2(T)]>0,$ for some $\boldsymbol{\theta}\in\Theta_0$, is inadmissible for estimating $\theta_2$ and is dominated by the estimator $\delta_{\psi_2^{*}}(\bold{X})=\psi_2^{*}(T)X_2$, where $\psi_2^{*}(T)=\max\big\{\frac{\alpha_1+\alpha_2+2}{\alpha_1+\alpha_2+1}\max\big\{1,\frac{1}{T}\big\},\psi_2(T)\big\}.$\vspace*{2mm}

The restricted MLE of $\theta_2$ is $\delta_{2,M}(\bold{X})=\max\{X_1,X_2\}$.
The restricted MLE $\delta_{2,M}(\bold{X})=\max\{X_1,X_2\}$ is inadmissible for estimating $\theta_2$ and is dominated by $\delta_{\psi_{2,M}^{*}}(\bold{X})=\frac{\alpha_1+\alpha_2+2}{\alpha_1+\alpha_2+1}\max\{X_1,X_2\}$.
\vspace*{2mm}

The unrestricted BSEE of $\theta_2$ is $\delta_{0,2}(\bold{X})=\frac{\alpha_2+2}{\alpha_2+1}X_2$. It follows that the unrestricted BSEE $\delta_{0,2}(\bold{X})=\frac{\alpha_2+2}{\alpha_2+1}X_2$ is inadmissible and is dominated by the estimator $\delta_{\psi_{0,2}^{*}}(\bold{X})=
\max\Big\{\frac{\alpha_2+2}{\alpha_2+1}X_2,\frac{\alpha_1+\alpha_2+2}{\alpha_1+\alpha_2+1}X_1\Big\}.$ Misra and Dhariyal (\citeyear{MR1326266}) derived similar results.
\\~\\\textbf{Example 3.2.3.} Let $X_1$ and $X_2$ be independent gamma random variables as defined in Example 3.1.3. 
Consider estimation of $\theta_2$, under the squared error loss function. From (\ref{eq:3.2.1}), we have 
$$\psi_{\lambda}(t)=\frac{\int_{0}^{\infty} s\,h_2(\frac{t}{\lambda}\vert s)f_2(s) ds}{\int_{0}^{\infty} s^2\,h_2(\frac{t}{\lambda}\vert s)f_2(s) ds}=\frac{\left(1+\frac{\lambda}{t}\right)}{\alpha_1+\alpha_2+1}, \;  t\in (0,\infty),\; \lambda\geq 1.$$
Clearly, for $t\in \Re_{++}$, $\underline{\psi}(t)=\frac{\left(1+\frac{1}{t}\right)}{1+\alpha_1+\alpha_2}\;\text{and}\; 
\overline{\psi}(t)=\infty.$\vspace*{2mm}

Let $\delta_{\psi_2}(\bold{X})=\psi_2(T) X_2$ be an arbitrary equivariant estimator of $\theta_2$. By Theorem 3.2.1, the equivariant estimator $\delta_{\psi_2}(\bold{X})=\psi_2(T) X_2$ is inadmissible and is dominated by the estimator $\delta_{\psi_2^{*}}(\bold{X})=\psi_2^{*}(T)X_2$, where
$\psi_2^{*}(T)=\max\big\{\psi_2(T),\frac{1+\frac{1}{T}}{\alpha_1+\alpha_2+1}\big\}$,  provided that $P_{\boldsymbol{\theta}}[\underline{\psi}(T)> \psi_2(T)]>0,$ for some $\boldsymbol{\theta}\in\Theta_0$.\vspace*{2mm}

Consequently, the unrestricted BSEE $\delta_{0,2}(\bold{X})=\frac{X_2}{\alpha_2+1}$ is inadmissible and an improved estimator is given by $\delta_{\psi_{0,2}^{*}}(\bold{X})=\max\big\{\frac{X_2}{\alpha_2+1},\,\frac{\left(X_1+X_2\right)}{\alpha_1+\alpha_2+1}\big\}$. Similar results are also obtained by Vijayasree et al. (\citeyear{MR1345425}).

\section{\textbf{Simulation Study}}\label{4}
\subsection{For Location Parameter $\theta_1$}
\noindent
\vspace*{2mm}

In Example 2.1.1, we have considered a bivariate normal distribution with unknown order restricted means (i.e., $\theta_1\leq \theta_2$), known variances ($\sigma_1^2$ and $\sigma_2^2$) and known correlation coefficient ($\rho$). For estimation of $\theta_1$ under the squared error loss function, we obtained estimators (given by (\ref{eq:2.1.3}) and (\ref{eq:2.1.4})) improving on the BLEE $X_1$. The improved estimators are same as the restricted maximum likelihood estimator. To further evaluate the performances of improved estimators, in this section, we compare the risk performances of the BLEE $X_1$ and the improved BLEE (as defined in (\ref{eq:2.1.3}) and (\ref{eq:2.1.4})), numerically, through Monte Carlo simulations. The simulated risks of the BLEE and the improved BLEE (restricted MLE) have been computed.

\begin{figure}
	\begin{subfigure}{.48\textwidth}
		\centering
		\includegraphics[width=68mm,scale=1.2]{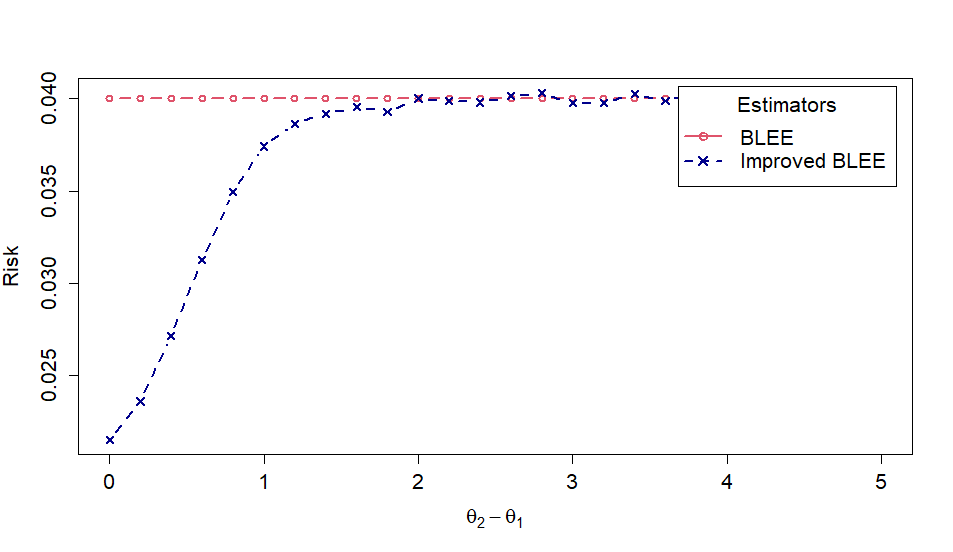} 
		\caption{$\sigma_1=0.2$, $\sigma_2=0.4$ and $\rho=-0.9$.} 
		\label{fig7:a} 
	\end{subfigure}
	\begin{subfigure}{.48\textwidth}
		\centering
		\includegraphics[width=68mm,scale=1.2]{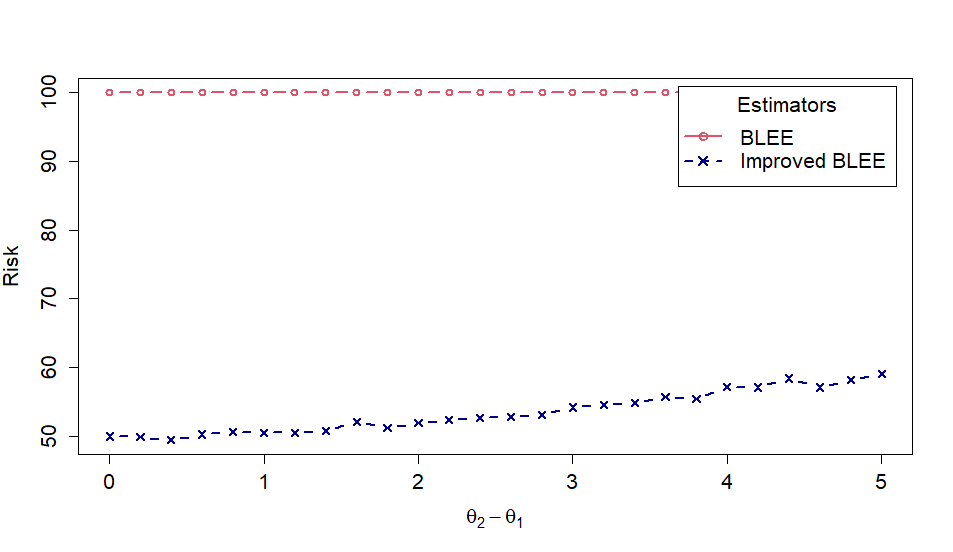} 
		
		\caption{$\sigma_1=10$, $\sigma_2=0.4$ and $\rho=-0.5$.} 
		\label{fig7:b} 
	\end{subfigure}
	\\	\begin{subfigure}{.48\textwidth}
		\centering
		\includegraphics[width=68mm,scale=1.2]{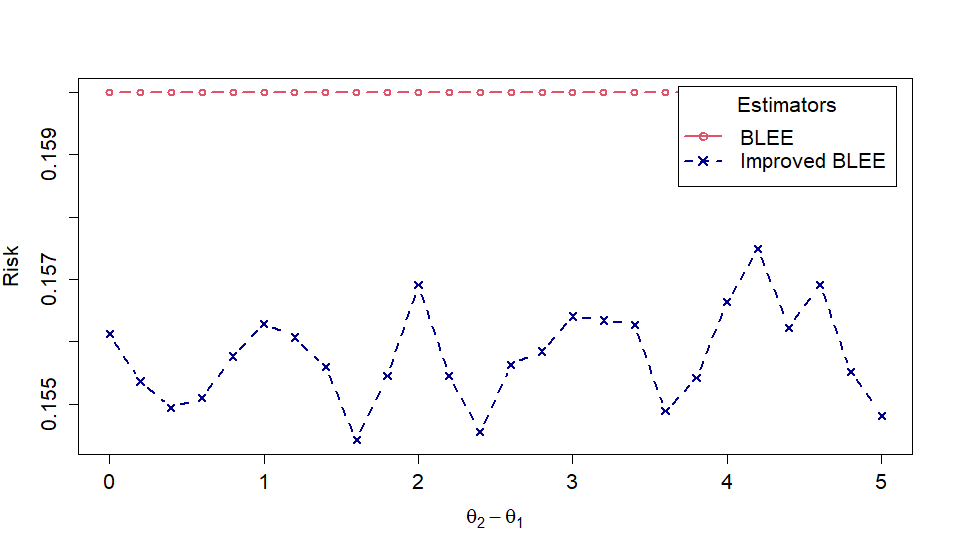} 
		\caption{$\sigma_1=0.4$, $\sigma_2=10$ and $\rho=-0.2$.} 
		\label{fig7:c} 
	\end{subfigure}
	\begin{subfigure}{.48\textwidth}
		\centering
		
		\includegraphics[width=68mm,scale=1.2]{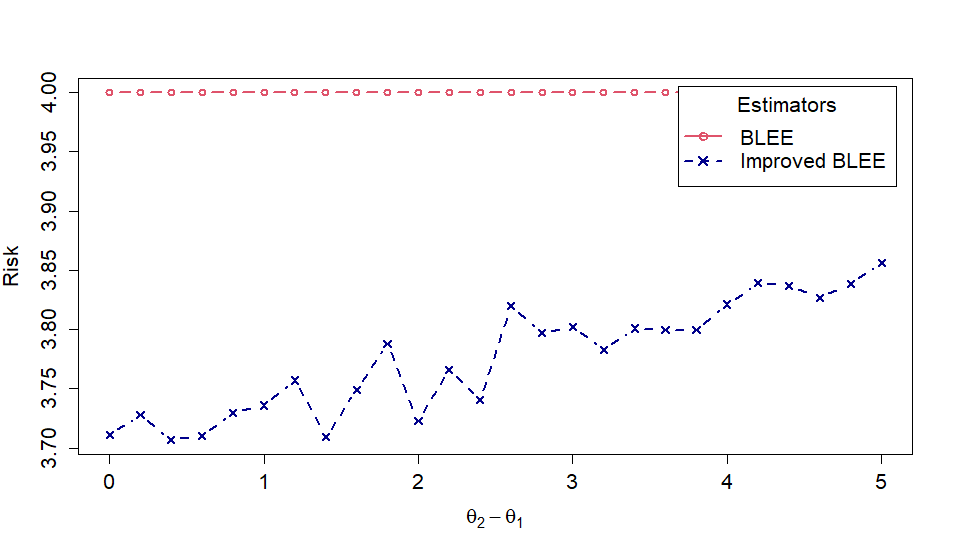} 
		
		\caption{$\sigma_1=2$, $\sigma_2=5$ and $\rho=0$.} 
		\label{fig7:d} 
	\end{subfigure}
	\\	\begin{subfigure}{.48\textwidth}
		\centering
		\includegraphics[width=68mm,scale=1.2]{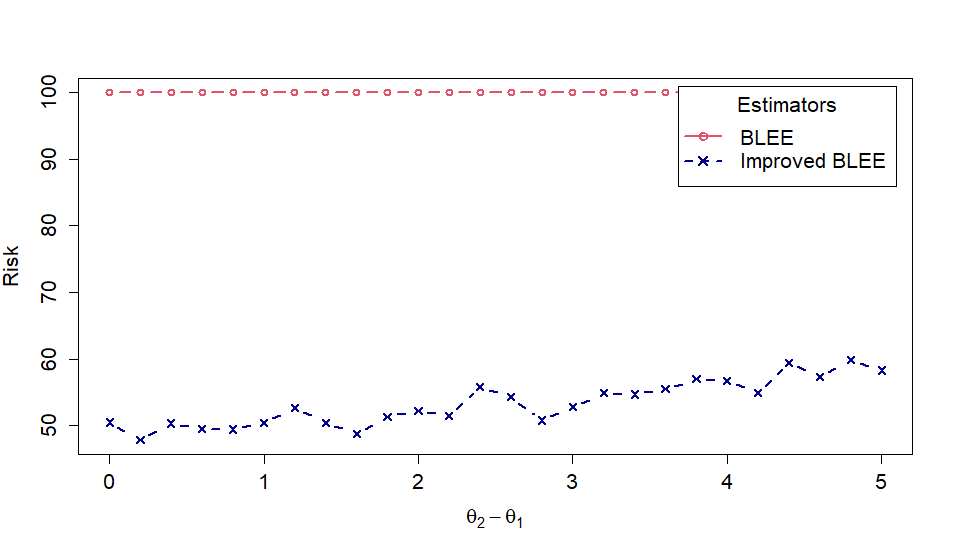} 
		\caption{ $\sigma_1=10$, $\sigma_2=0.4$ and $\rho=0$.} 
		\label{fig7:e}  
	\end{subfigure}
	\begin{subfigure}{.48\textwidth}
		\centering
		
		\includegraphics[width=68mm,scale=1.2]{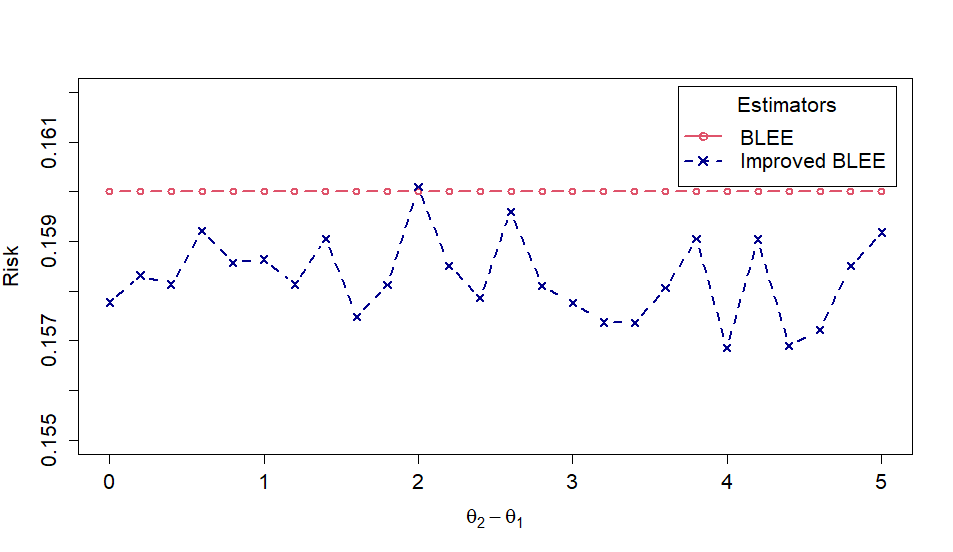} 
		\caption{$\sigma_1=0.4$, $\sigma_2=10$ and $\rho=0.2$.} 
		\label{fig7:f} 
	\end{subfigure}
	\\	\begin{subfigure}{.48\textwidth}
		\centering
		\includegraphics[width=68mm,scale=1.2]{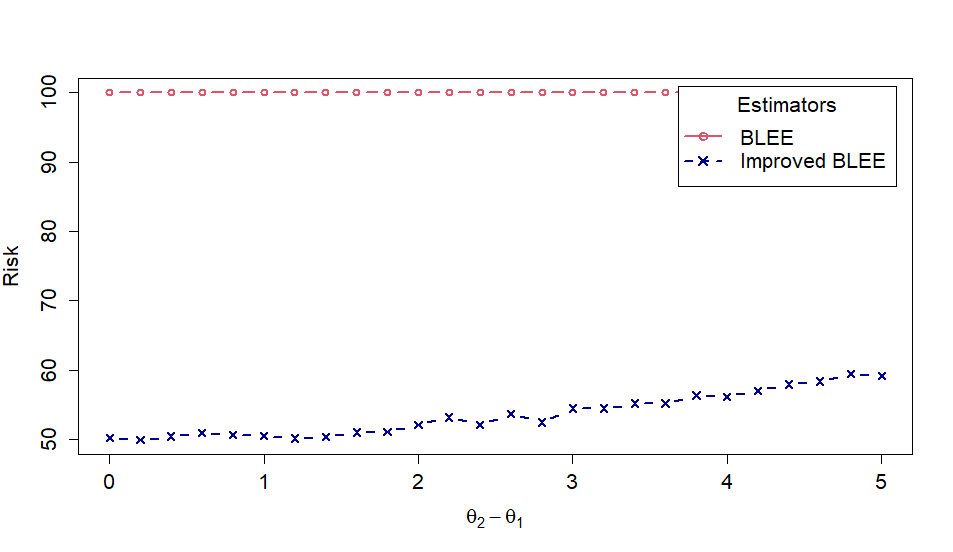} 
		\caption{ $\sigma_1=10$, $\sigma_2=0.4$ and $\rho=0.5$.} 
		\label{fig7:g}  
	\end{subfigure}
	\begin{subfigure}{.48\textwidth}
		\centering
		
		\includegraphics[width=68mm,scale=1.2]{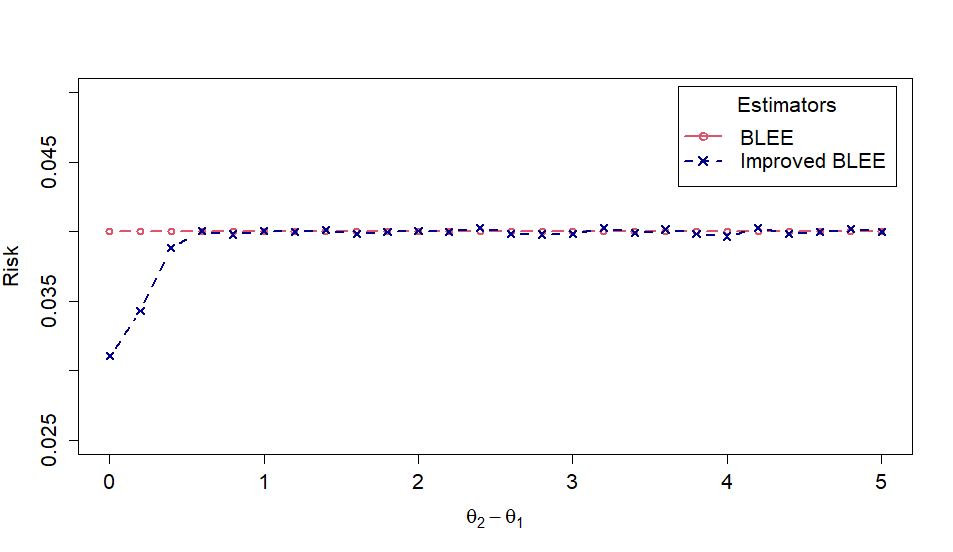} 
		\caption{$\sigma_1=0.2$, $\sigma_2=0.4$ and $\rho=0.9$.} 
		\label{fig7:h} 
	\end{subfigure}
	
	\caption{Risk plots of Estimators of smaller location parameter $\theta_1$ against values of $\theta_2-\theta_1$: Bivariate Normal distribution.}
	\label{fig7}
\end{figure}

Note that the BLEE has constant rsik $\sigma_1^2$. For simulations, 10000 random samples of size 1 were generated from the relevant bivariate normal distribution. The simulated values of the risks of the improved BLEE against the constant risk $\sigma_1^2$ of the BLEE are plotted in Figure \ref{fig7}. The following observations are evident from Figure \ref{fig7}:
\\(i) The improved BLEE performs better than the BLEE, which is in conformity with the theoretical findings of Example 2.1.1.
\\(ii) The performance of the improved BLEE is significantly better than the BLEE when $\theta_1$ and $\theta_2$ ($\theta_1\leq \theta_2$) are close.
\\(iii) For larger values of $\sigma_1$ and smaller values of $\sigma_2$ the improved BLEE performs significantly better than the BLEE even when $\theta_1$ and $\theta_2$ ($\theta_1\leq \theta_2$) are not close. Similarly, for smaller values of $\sigma_1$, larger values of $\sigma_2$ and non-positive $\rho$, the improved BLEE outperforms the BLEE.

\subsection{For Scale Parameter $\theta_1$}
\noindent
\vspace*{2mm}

Example 3.1.3 concerns estimation of the order restricted scale parameters of two gamma distributions under the squared error loss function. The shape parameter $\alpha_1>0$ and $\alpha_2>0$ are assumed to be known. The risk performances of the restricted MLE $\min\!\big\{\frac{X_1}{\alpha_1},\frac{X_1+X_2}{\alpha_1+\alpha_2}\}$, the BSEE $\frac{X_1}{\alpha_1+1}$, the improved restricted MLE $\min\!\big\{\frac{X_1}{\alpha_1},\frac{X_1+X_2}{\alpha_1+\alpha_2+1}\}$ and the improved BSEE $\min\!\big\{\frac{X_1}{\alpha_1+1},\frac{X_1+X_2}{\alpha_1+\alpha_2+1}\}$, under the squared error loss function, are compared through a simulation study.

\FloatBarrier
\begin{figure}[h!]
	
	\begin{subfigure}{.48\textwidth}
		\centering
		\includegraphics[width=68mm,scale=1.2]{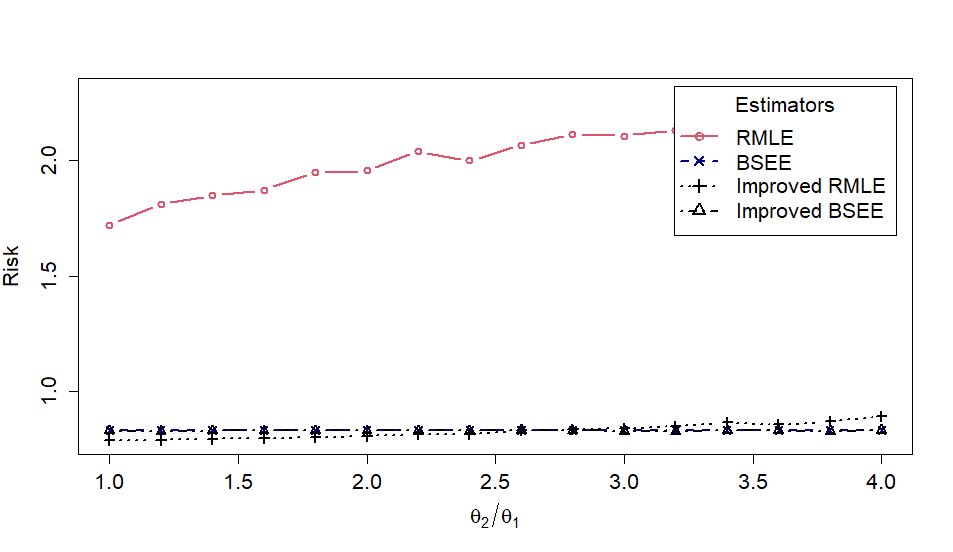} 
		\caption{$\alpha_1=0.2$ and $\alpha_2=0.2$.} 
		\label{fig8:a} 
	\end{subfigure}
	\begin{subfigure}{.48\textwidth}
		\centering
		\includegraphics[width=68mm,scale=1.2]{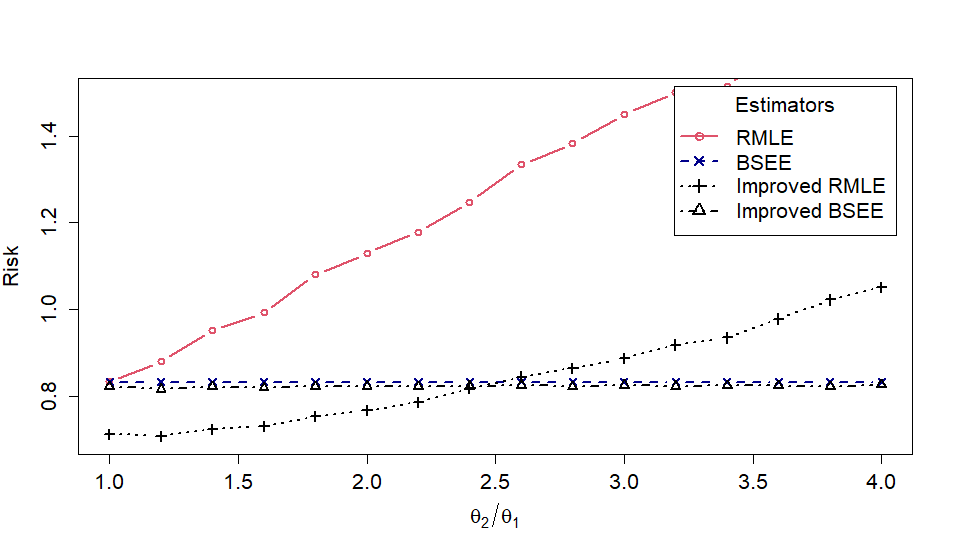} 
		
		\caption{$\alpha_1=0.2$ and $\alpha_2=1$.} 
		\label{fig8:b} 
	\end{subfigure}
	\\	\begin{subfigure}{.48\textwidth}
		\centering
		\includegraphics[width=68mm,scale=1.2]{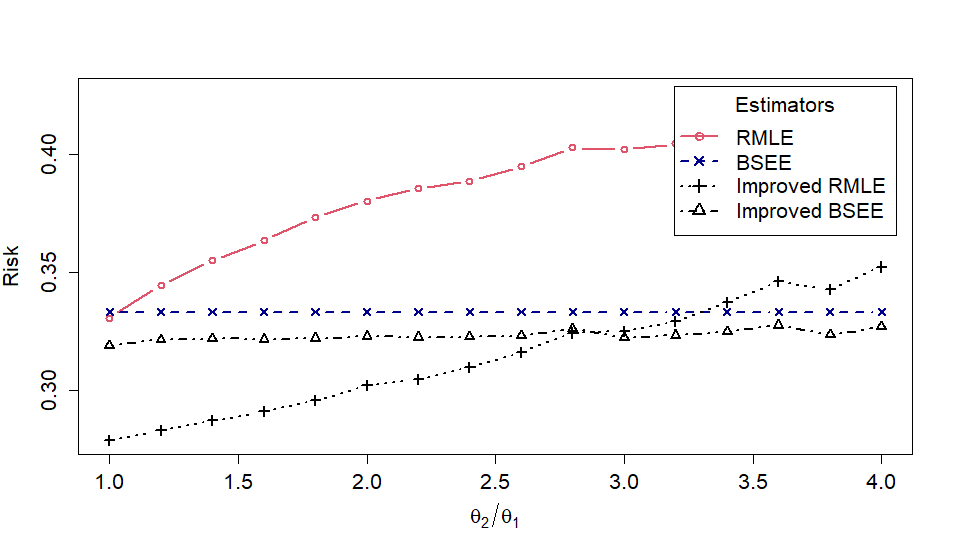} 
		\caption{$\alpha_1=2$ and $\alpha_2=1$.} 
		\label{fig8:c} 
	\end{subfigure}
	\begin{subfigure}{.48\textwidth}
		\centering
		
		\includegraphics[width=68mm,scale=1.2]{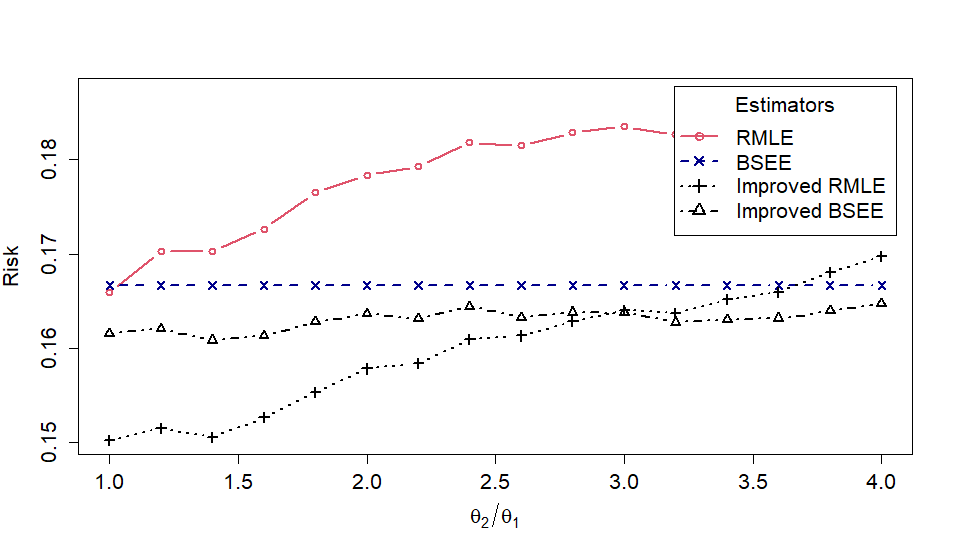} 
		
		\caption{$\alpha_1=5$ and $\alpha_2=1$.} 
		\label{fig8:d}
	\end{subfigure}
	\\	\begin{subfigure}{.48\textwidth}
		\centering
		\includegraphics[width=68mm,scale=1.2]{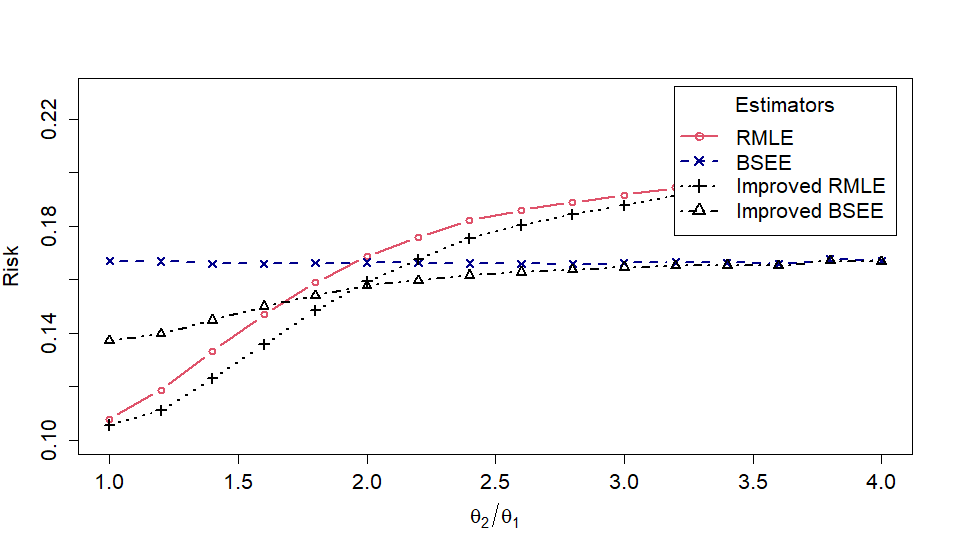} 
		\caption{$\alpha_1=5$ and $\alpha_2=10$.} 
		\label{fig8:e} 
	\end{subfigure}
	\begin{subfigure}{.48\textwidth}
		\centering
		
		\includegraphics[width=68mm,scale=1.2]{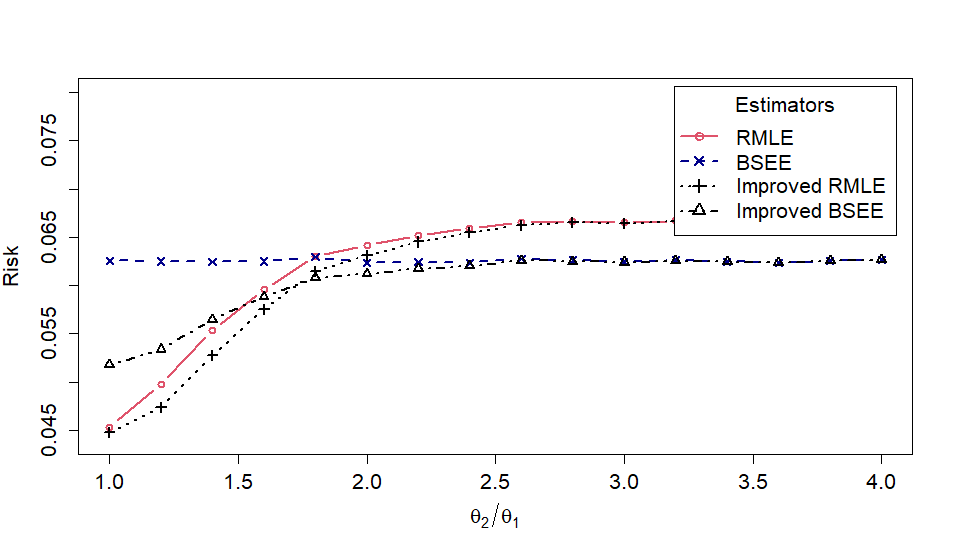} 
		\caption{$\alpha_1=15$ and $\alpha_2=15$.} 
		\label{fig8:f} 
	\end{subfigure}
	
	\caption{Risk plots of Estimators of smaller scale parameter $\theta_1$ against values of $\theta_2/\theta_1$: Independent Gamma distributions.}
	\label{fig8} 
\end{figure}

Note that the BSEE $\frac{X_1}{\alpha_1+1}$ has the constant risk $\frac{1}{\alpha_1+1}$. The simulated risk functions of the restricted MLE (RMLE), the BSEE, the improved restricted MLE and the improved BSEE are computed using 10000 random samples of size 1, generated from relevant gamma distributions. The simulated risks of various estimators are plotted in Figure \ref{fig8}. The following observations are evident from Figure \ref{fig8}:
\\	(i) The improved BSEE and the improved restricted MLE always perform better than the BSEE and the restricted MLE, respectively, which is in conformity with theoretical findings of Example 3.1.3.
\\(ii) We also observe that the improved restricted MLE performs better than the improved BSEE, for small values of $\frac{\theta_2}{\theta_1}$, and the improved BSEE dominates the improved restricted MLE for the moderate and large values of $\frac{\theta_2}{\theta_1}$. 
\\(iii) For smaller values of shape parameters $\alpha_1$ and $\alpha_2$ the performance of the restricted MLE is worse than other estimators.

\section{\textbf{Real Life Data Analysis}} \label{5}

We first consider modelling the time taken by male professional sprinters to complete 200m and 100m sprints. To build the model, we studied the data of winners of the men's Olympic sprint finals from 1900 to 2022 in 200m and 100m sprints, as available on websites statista1\footnote{https://www.statista.com/statistics/1091753/olympics-200m-gold-medal-times-since-1900/} and statista2\footnote{https://www.statista.com/statistics/1090316/olympics-100m-gold-medal-times-since-1896/}. To test the bivariate normality of this data set, we applied the Shapiro-Wilk and Henze-Zirkler multivariate normality tests and observed p-values of $0.358$ and $0.216$, respectively. This suggests that there is no significant departure from the normality. Let $\theta_1$ and $\theta_2$ be the means of athletes' average speeds in 200m and 100m sprints, respectively. We then used the paired t-test for the hypotheses $H_0: \theta_1 \geq \theta_2\; \text{vs }H_1: \theta_1< \theta_2$, and observed the p-value was $0.11$. In view of this, it is reasonable to assume that $\theta_1\leq \theta_2$. This is not surprising as one would anticipate a sprinter to slow down after reaching top speed somewhere in the first 100 meters. Also, Usain Bolt's 200m and 100m world record times of 19.19 seconds and 9.58 seconds work out to his average speeds of 10.42 m/s and 10.44 m/s, respectively, in 200m and 100m sprints. The above analysis suggests that for a data on average speeds of professional sprinters in finishing 200m and 100m sprints it is reasonable to assume bivariate normality with order restricted means.\vspace*{2mm}

Now, to illustrate an application of the findings of our paper, another data set of male professional sprinters named "UK Sprinters data" is considered. The details of this data set are available in the book "A Handbook of Small Data Sets" by D.J. Hand et al. (\citeyear{hand}) (p. 53). This data is presented in Table \ref{table:1}. The sample means, sample variances and sample correlation coefficient of this data are $9.617,\;9.488,\; 0.032,\;0.063\;\text{and}\;0.848$, respectively. In view of the above modelling, we assume that this data follows a bivariate normal distribution with order restricted means. Consequently, the sample means of average speed (m/s) of 200m and 100m sprints also obey bivariate normality with order restricted means. Let $X_1$ and $X_2$ be random variables that represent the sample means of average speed (m/s) of 200m and 100m sprints, respectively. Then, $(X_1,X_2)$ follows a bivariate normal distribution with means $\theta_1$ and $\theta_2$, variances $\sigma_1^2$ and $\sigma_2^2$, and correlation coefficient $\rho$, where $\theta_1\leq \theta_2$. We use the sample variances and sample correlation coefficient of the data of Table \ref{table:1} to obtain the plug-in values for $\sigma_1^2$, $\sigma_2^2$ and $\rho$, as $\sigma_1^2=\frac{0.032}{11}=0.0029$, $\sigma_2^2=\frac{0.063}{11}=0.0057$ and $\rho=0.848$. In the absence of prior information $\theta_1\leq \theta_2$, one would like to use $x_1=9.617$ and $x_2=9.488$ as estimates of $\theta_1$ and $\theta_2$, respectively. Now using the information that $\theta_1 \leq \theta_2$, from Example 2.1.1 and Example 2.2.1, we can obtain improvements over unrestricted estimators $x_1=9.617$ and $x_2=9.488$ for estimating $\theta_1$ and $\theta_2$, receptively. Note that here $\rho \sigma_2 >\sigma_1$ and $\rho \sigma_1 <\sigma_2$. By Example 2.1.1 and Example 2.2.1, the estimates
$$ \max\Big\{x_1, \frac{\sigma_2(\sigma_2-\rho\sigma_1)x_1+\sigma_1(\sigma_1-\rho\sigma_2)x_2}{\sigma_1^2+\sigma_2^2-2\rho \sigma_1 \sigma_2}\Big\}=9.66$$
and
$$ \max\Big\{x_2, \frac{\sigma_2(\sigma_2-\rho\sigma_1)x_1+\sigma_1(\sigma_1-\rho\sigma_2)x_2}{\sigma_1^2+\sigma_2^2-2\rho \sigma_1 \sigma_2}\Big\}=9.66$$
are better estimates of $\theta_1$ and $\theta_2$, respectively. Based on the above analysis, we infer that, rather than using 9.617 and 9.488 as estimates of $\theta_1$ and $\theta_2$, respectively, 9.66 be taken as a common estimated value for $\theta_1$ and $\theta_2$.

\begin{table}[h!]
	\centering
	\begin{tabular}{|c | c|  c | c | c|} 
		\hline
		Athlete & Time in 200m & Time in 100m & 200m Average & 100m Average \\ 
		Name &  (in seconds) & (in seconds) & speed (m/s)  & speed (m/s) \\
		[0.5ex] 
		\hline\hline
		L Christie & 20.09 & 9.97 & 9.955 & 10.030\\ 
		J Regis & 20.32 & 10.31 & 9.842 & 9.699\\
		M Rossowess & 20.51 & 10.40 & 9.751 & 9.615\\
		A Carrott & 20.76 & 10.56 & 9.633 & 9.469\\
		T Bennett & 20.90 & 10.92 & 9.569 & 9.157\\
		A Mafe & 20.94 & 10.64 & 9.551 & 9.398\\ 
		D Reid & 21.00 & 10.54 & 9.523 & 9.487\\
		P Snoddy & 21.14 & 10.85 & 9.461 & 9.216\\
		L Stapleton & 21.17 & 10.71 & 9.447 & 9.337\\
		C Jackson & 21.19 & 10.56 & 9.438 & 9.469\\ 
		\hline\hline
		Mean & &  &$x_1$=9.617 &$x_2$= 9.488\\
		\hline
	\end{tabular}
	\caption{UK Sprinters Data Set}
	\label{table:1}
\end{table}

\section{\textbf{Concluding Remarks}}
For the problem of estimating order restricted location/scale parameters of two distributions, we unify various results proved in the literature for specific probability distributions, having independent marginal distributions, and specific loss functions. We unify these studies by considering a general bivariate probability location/scale model (possibly having dependent marginals) and a quite general loss function. We derive sufficient conditions for the inadmissibility of an arbitrary location/scale equivariant estimator. For estimators satisfying these sufficient conditions, we also obtain dominating estimators. It may be possible to extend this study to estimate $k\, (\geq 2)$ order restricted parameters as in Vijayasree et al. (\citeyear{MR1345425}). This will be considered in our future research.

\section*{\textbf{Funding}}

This work was supported by the [Council of Scientific and Industrial Research (CSIR)] under Grant [number 09/092(0986)/2018].

\bibliographystyle{apalike}
\bibliography{Paper2}

\begin{thebibliography}{}

\bibitem[Barlow et~al., 1972]{MR0326887}
Barlow, R.~E., Bartholomew, D.~J., Bremner, J.~M., and Brunk, H.~D. (1972).
\newblock {\em Statistical inference under order restrictions. {T}he theory and
  application of isotonic regression}.
\newblock John Wiley \& Sons.

\bibitem[Garg and Misra, 2021]{garg2021componentwise}
Garg, N. and Misra, N. (2021).
\newblock Componentwise equivariant estimation of order restricted location and
  scale parameters in bivariate models: A unified study.

\bibitem[Garren, 2000]{MR1802627}
Garren, S.~T. (2000).
\newblock On the improved estimation of location parameters subject to order
  restrictions in location-scale families.
\newblock {\em Sankhy\={a} Ser. B}, 62(2):189--201.

\bibitem[Hand et~al., 1993]{hand}
Hand, D.~J., Daly, F., McConway, K., Lunn, D., and Ostrowski, E. (1993).
\newblock {\em A handbook of small data sets}.
\newblock cRc Press.

\bibitem[Kaur and Singh, 1991]{MR1128873}
Kaur, A. and Singh, H. (1991).
\newblock On the estimation of ordered means of two exponential populations.
\newblock {\em Ann. Inst. Statist. Math.}, 43(2):347--356.

\bibitem[Kotz et~al., 2000]{MR1788152}
Kotz, S., Balakrishnan, N., and Johnson, N.~L. (2000).
\newblock {\em Continuous multivariate distributions. {V}ol. 1}.
\newblock Wiley Series in Probability and Statistics: Applied Probability and
  Statistics. Second edition.

\bibitem[Kubokawa and Saleh, 1994]{MR1370413}
Kubokawa, T. and Saleh, A. K. M.~E. (1994).
\newblock Estimation of location and scale parameters under order restrictions.
\newblock {\em J. Statist. Res.}, 28(1-2):41--51.

\bibitem[Kumar and Sharma, 1988]{MR981031}
Kumar, S. and Sharma, D. (1988).
\newblock Simultaneous estimation of ordered parameters.
\newblock {\em Comm. Statist. Theory Methods}, 17(12):4315--4336.

\bibitem[Misra and Dhariyal, 1995]{MR1326266}
Misra, N. and Dhariyal, I.~D. (1995).
\newblock Some inadmissibility results for estimating ordered uniform scale
  parameters.
\newblock {\em Comm. Statist. Theory Methods}, 24(3):675--685.

\bibitem[Patra and Kumar, 2017]{Patra}
Patra, L.~K. and Kumar, S. (2017).
\newblock Estimating ordered means of a bivariate normal distribution.
\newblock {\em American Journal of Mathematical and Management Sciences},
  36(2):118--136.

\bibitem[Robertson et~al., 1988]{MR961262}
Robertson, T., Wright, F.~T., and Dykstra, R.~L. (1988).
\newblock {\em Order restricted statistical inference}.
\newblock John Wiley \& Sons.

\bibitem[Stein, 1964]{MR171344}
Stein, C. (1964).
\newblock Inadmissibility of the usual estimator for the variance of a normal
  distribution with unknown mean.
\newblock {\em Ann. Inst. Statist. Math.}, 16:155--160.

\bibitem[Taketomi et~al., 2021]{axioms10040267}
Taketomi, N., Konno, Y., Chang, Y.-T., and Emura, T. (2021).
\newblock A meta-analysis for simultaneously estimating individual means with
  shrinkage, isotonic regression and pretests.
\newblock {\em Axioms}, 10(4).

\bibitem[Van~Eeden, 2006]{MR2265239}
Van~Eeden, C. (2006).
\newblock {\em Restricted parameter space estimation problems}, volume 188 of
  {\em Lecture Notes in Statistics}.
\newblock Springer, New York.

\bibitem[Vijayasree et~al., 1995]{MR1345425}
Vijayasree, G., Misra, N., and Singh, H. (1995).
\newblock Componentwise estimation of ordered parameters of {$k\ (\geq 2)$}
  exponential populations.
\newblock {\em Ann. Inst. Statist. Math.}, 47(2):287--307.

\end{thebibliography}

\end{document}